\documentclass[12pt,leqno]{article}
\usepackage{amsmath,amsthm,amscd,amsfonts,amssymb}

\newcounter{num}
\newtheorem{theorem}{Theorem}[section]

\newtheorem{corollary}[theorem]{Corollary}
\newtheorem{proposition}[theorem]{Proposition}
\newtheorem{lemma}[theorem]{Lemma}

\newtheorem{remark}[theorem]{Remark}

\begin{document}

\baselineskip 22pt plus 2pt

\newcommand{\be}        {\begin{eqnarray}}
\newcommand{\ee}        {\end{eqnarray}}
\newcommand{\pl}{\partial}
\newcommand{\sbs}{\subset}
\newcommand{\vr}{\varphi}

\newcommand{\AAA}       {{\Bbb A}}     % NB: A,E,G,L,S,T are triple letters
\newcommand{\aaa}       {{\LBbb A}}
\newcommand{\BB}        {{\Bbb B}}
\newcommand{\bb}        {{\LBbb B}}
\newcommand{\CC}        {{\Bbb C}}
\newcommand{\cc}        {{\LBbb C}}
\newcommand{\DD}        {{\Bbb D}}
\newcommand{\dd}        {{\LBbb D}}
\newcommand{\EEE}       {{\Bbb E}}
\newcommand{\eee}       {{\LBbb E}}
\newcommand{\FF}        {{\Bbb F}}
\newcommand{\ff}        {{\LBbb F}}
\newcommand{\GGG}       {{\Bbb G}}
\newcommand{\HH}        {{\Bbb H}}
\newcommand{\hh}        {{\LBbb H}}
\newcommand{\II}        {{\Bbb I}}
\newcommand{\ii}        {{\LBbb I}}
\newcommand{\JJ}        {{\Bbb J}}
\newcommand{\jj}        {{\LBbb J}}
\newcommand{\KK}        {{\Bbb K}}
\newcommand{\kk}        {{\LBbb k}}
\newcommand{\LLL}       {{\Bbb L}}
\newcommand{\MM}        {{\Bbb M}}
\newcommand{\mm}        {{\LBbb M}}
\newcommand{\NN}        {{\Bbb N}}
\newcommand{\nn}        {{\LBbb N}}
\newcommand{\OO}        {{\Bbb O}}
\newcommand{\oo}        {{\LBbb O}}
\newcommand{\PP}        {{\Bbb P}}
\newcommand{\pp}        {{\LBbb P}}
\newcommand{\QQ}        {{\Bbb Q}}
\newcommand{\qq}        {{\LBbb Q}}
\newcommand{\RR}        {{\Bbb R}}
\newcommand{\rr}        {{\LBbb R}}
\newcommand{\SSS}       {{\Bbb S}}
\newcommand{\sss}       {{\LBbb S}}
\newcommand{\TTT}       {{\Bbb T}}
\newcommand{\ttt}       {{\LBbb t}}
\newcommand{\UU}        {{\Bbb U}}
\newcommand{\VV}        {{\Bbb V}}
\newcommand{\vv}        {{\LBbb V}}
\newcommand{\WW}        {{\Bbb W}}
\newcommand{\ww}        {{\LBbb W}}
\newcommand{\XX}        {{\Bbb X}}
\newcommand{\xx}        {{\LBbb X}}
\newcommand{\YY}        {{\Bbb Y}}
\newcommand{\yy}        {{\LBbb Y}}
\newcommand{\ZZ}        {{\Bbb Z}}
\newcommand{\zz}        {{\LBbb Z}}

\newcommand{\calA}      {{\cal A}}
\newcommand{\calB}      {{\cal B}}
\newcommand{\calC}      {{\cal C}}
\newcommand{\calD}      {{\cal D}}
\newcommand{\calE}      {{\cal E}}
\newcommand{\calF}      {{\cal F}}
\newcommand{\calG}      {{\cal G}}
\newcommand{\calH}      {{\cal H}}
\newcommand{\calI}      {{\cal I}}
\newcommand{\calJ}      {{\cal J}}
\newcommand{\calK}      {{\cal K}}
\newcommand{\calL}      {{\cal L}}
\newcommand{\calM}      {{\cal M}}
\newcommand{\calN}      {{\cal N}}
\newcommand{\calO}      {{\cal O}}
\newcommand{\calP}      {{\cal P}}
\newcommand{\calQ}      {{\cal Q}}
\newcommand{\calR}      {{\cal R}}
\newcommand{\calS}      {{\cal S}}
\newcommand{\calT}      {{\cal T}}
\newcommand{\calU}      {{\cal U}}
\newcommand{\calV}      {{\cal V}}
\newcommand{\calW}      {{\cal W}}
\newcommand{\calX}      {{\cal X}}
\newcommand{\calY}      {{\cal Y}}
\newcommand{\calZ}      {{\cal Z}}

\begin{titlepage}
\vspace*{1cm}

\begin{center}
{\Large\bf Isometries of Hilbert $C^*$-modules} \vspace{2cm}

{\large Baruch Solel}\footnote{Supported by Technion V.P.R. Fund--Steigman Research Fund,
Technion V.P.R. Fund--Fund for the Promotion of Sponsored Reserach and the Fund for the
Promotion of Research at the Technion.}

Department of Mathematics\\[-0.2cm]
Technion -- Israel Institute of Technology\\[-0.2cm]
Haifa 32000\\[-0.2cm]
Israel

\end{center}
\vspace{2cm}

%                       Abstract

\begin{abstract}

Let $X$ and $Y$ be right, full, Hilbert $C^*$-modules over the
algebras $A$ and $B$ respectively and let $T:X\to Y$ be a linear
surjective isometry . Then $T$ can be extended to an isometry of
the linking algebras. $T$ then is a sum of two maps: a (bi-)module
map (which is completely isometric and preserves the inner
product) and a  map that reverses the (bi-)module actions. If $A$
(or $B$) is a factor von Neumann algebra then every isometry
$T:X\to Y$ is either a (bi-)module map or reverses the (bi-)module
actions.
\end{abstract}

\end{titlepage}

%                           Introduction

\section{Introduction}

Given a right Hilbert $C^*$-module $X$ over a $C^*$-algebra $A$ it
is a module over $A$ and has an $A$-valued inner product.  One
then defines the norm of $X$ using the inner product and it makes
$X$ a Banach space. It is known that once the module structure and
Banach space structure are given (for a $C^*$-module $X$) the
$A$-valued inner product is uniquely defined. This was proved by
Lance in \cite[Theorem]{L1} and, independently by Blecher in
\cite[Theorems~3.1 and 3.2]{B1}. In fact, as Blecher showed, the
inner product can be recovered from the module and Banach space
structures. This result of Lance and Blecher can be stated as
follows.

%                     Theorem 1.1

\begin{theorem}[\cite{B2}, \cite{L1}]
Let $X_1$ and $X_2$ be right Hilbert $C^*$-modules over a
$C^*$-algebra $A$ and let $S:X_1\to X_2$ be a surjective isometry
which is an $A$-module map. Then $S$ preserves the inner product,
i.e. $\left<Sx,Sy \right>_2=\left<x,y \right> _1$ (where
$\left<\cdot,\cdot \right>_j$ is the inner product in $X_j$).
Moreover, the inner product of a right Hilbert $C^*$-module $X$
over $A$ can be recovered from the norm and the module structure
by
$$
\left<x,x\right>=\sup\{r(x)^*r(x):r \ {\rm is \ an} \ A{\rm -module \ map:}
X\to A, \|r\|\le 1 \}
$$
and
$$
\left<x,y \right>=\frac 14 \sum^3_{k=0} i^k \left<x+i^ky, \ x+i^ky \right>
\ \ \ (i=\sqrt{-1}).
$$

\end{theorem}

Another proof can be found in \cite[Theorem~5]{F}. One can modify
the first part of the theorem for the case where $X_1$ is a
$C^*$-module over $A$ and $X_2$ is a $C^*$-module over $B$ and $S$
is a module map in the sense that there is a $*$-isomorphism $\alpha: A\to B$
such that $S(xa)=(Sx)\alpha(a)$. In this case $S$ satisfies
$\left<Sx,Sy\right>_2=\alpha(\left<x,y \right>_1)$. (See
\cite[Lemma~5.10]{MS}).

In the present paper we study to what extent it is possible to
recover the $C^*$-module structure from the Banach space structure
alone. In other words, given an isometry $T$ (linear and
surjective) of a $C^*$-module $X$ over $A$ onto a $C^*$-module $Y$
over $B$, is it a module map? i.e. can we find a $^*$-isomorphism
$\alpha$ of $A$ onto $B$ with $T(xa)=T(x)\alpha(a)$? If we can,
then, as we mentioned above, we have $\left<Tx,Ty\right>
=\alpha(\left<x,y \right>)$. Also, denoting by $y\otimes z^*$ (for
$y,z\in X$) the operator on $X$ defined by
$(y\otimes z^*)(x)=y\left<z,x\right>$, we have
\be
T\bigl((y\otimes z^*)(x)\bigr)=T\bigl(y\left<z,x\right> \bigr)=
T(y)\alpha\bigl(\left<z,x\right>\bigr)\nonumber\\
=T(y)\left<T(z),T(x)\right>=\bigl(T(y)\otimes T(z)^*\bigr)(T(x)).\nonumber
\ee
The operators of the form $y\otimes z^*, \ y,z\in X$, generate a $C^*$-algebra
denoted by $\KK(X)$. It is possible to show that, for $T$ as above,
$\beta(y\otimes z^*)=Ty\otimes (Tz)^*$ defines a $^*$-isomorphism
$\beta:\KK(X)\to \KK(Y)$. The computation above shows that $T(Kx)=\beta(K)Tx$
\, for \, $x\in X, \ K\in\KK(X)$. The relationship we have now between $T,\alpha$
on $\beta$ can be summarized by considering the $^*$-algebra ${\cal L}(X)$ defined by
$$
{\cal L}(X)=
\begin{pmatrix}
\KK(X) & X\\
\bar X & A
\end{pmatrix}
$$
(where $\bar X$ and the product and involution on ${\cal L}(X)$ will be defined
shortly) and noting that the map $\psi:{\cal L}(X)\to {\cal L}(Y)$ defined by
$$
\psi
\begin{pmatrix}
K & x\\
\bar y & a
\end{pmatrix}=
\begin{pmatrix}
\beta(K) & Tx\\
\overline{Ty} & \alpha(a)
\end{pmatrix}
$$
is a $^*$-isomorphism.

Hence to say that $T$ preserves the $C^*$-module structure amounts to saying that
$T$ can be extended to a $^*$-isomorphism of ${\cal L}(X)$ onto ${\cal L}(Y)$.

By considering the transpose map that maps the Hilbert column space $H^c$ (a right
$C^*$-module over $\CC$, isometric to a Hilbert space $H$) onto the Hilbert row
space $H^r$ (a right Hilbert $C^*$-module over $K(H)$) it is clear that we don't
always have such a $^*$-isomorphism.

Our main result, Theorem~3.2, shows that, if $X$ and $Y$ are full, $T$ can always
be extended to an isometry of ${\cal L}(X)$ onto ${\cal L}(Y)$.

The celebrated result of Kadison~\cite[Theorem~7]{K} states that every unital isometry
of unital $C^*$-algebras is a selfadjoint Jordan map. For von Neumann algebras we can,
in fact, decompose the algebras as a direct sum of two summands. On one summand the map
is a $^*$-isomorphism and on the other it is a $^*$-antiisomorphism
(\cite[Theorem~10]{K}).
A similar result was proved also for isometries of some nonselfadjoint operator
algebras (\cite{S}). For an isometry $T$ of selfdual $C^*$-modules over von Neumann
algebras we find that $T$ can be written as a sum of an isometry which is a module
map (and preserves the inner product) and an isometry that is, in some sense, an
anti-module-map. (For a precise statement see Corollary~2.25).
The case of (not necessarily selfdual) Hilbert $C^*$-modules over  general
$C^*$-algebras is similar except that the decomposition of $X$ is done by a projection
in the enveloping von Neumann algebra of ${\cal L}(X)$ (Theorem~3.2).

As a corollary we show that, if we assume that the isometry $T$ is
in fact a 2-isometry (i.e., the map $I\otimes T:M_2\otimes X\to
M_2\otimes Y$ is an isometry), then $T$ preserves the $C^*$-module
structure (corollary~3.3). In particular, a 2-isometry of Hilbert
$C^*$-modules is necessarily a complete isometry.

After this work was completed it was pointed to us by D.~Blecher that M.~Hamana had
previously proved it~\cite{Ha} using different methods.

Also we show that, for a given Hilbert space $H, \ H^c$ and $H^r$
are the only Hilbert $C^*$-modules that are isometric to $H$
(Corollary~3.6).

Now we turn to set some notation and recall the
definitions that we need.
\medskip

%                        Definitions

\noindent
{\bf Definitions}\\
{\it
(1)~~A right \underline{ pre-Hilbert $C^*$-module} over a $C^*$-algebra $A$ is a
right-module $X$ equipped with a map \, $\left<\cdot,\cdot\right>:X\times X\to A$
satisfying:
\begin{description}
\item(i)~~~$\left< x,x \right>\ge 0, \ x\in X$ \, and \, $\left< x,x \right>= 0$ \,
only if \, $x=0$.

\item(ii)~~$\left< x,y \right>^*=\left< y,x \right> \ \ y,x\in X$.

\item(iii)~$y\longmapsto\left< x,y \right>$ is a linear map for all $x\in X$.

\item(iv)~~$\left< x,ya \right>=\left< x,y \right>a, \ \ x,y\in X, \ a\in A$.
\end{description}

\noindent
(2)~~The norm on a pre-Hilbert $C^*$-module $X$ over $A$ is defined by
$\|x\|=\|\left< x,x \right>\|^{\frac 12}$. If \ $X$ \ is complete with respect to this
norm then \ $X$\ is said to be a (right) \linebreak \underline{Hilbert  $C^*$-module}
over $A$.
\medskip

\noindent (3)~~A Hilbert  $C^*$-module $X$ over $A$ is said to be
\underline{full} if $A=\overline{{\rm span}}\{\left< x,y
\right>:x,y\in X\}$.  }
\medskip
\vspace{-0.5cm}

One can define left Hilbert $C^*$-module similarly. $X$ is then a
left $A$-module and the inner product is assumed to be linear in
the first entry. Also $\left< ax,y \right>=a\left< x,y \right>$.

Given a right Hilbert $C^*$-module $X$ over $A$ we define $\bar
X$, the {\it conjugate module}, as follows. As a set we write
$\bar X=\{\bar x:x\in X\}$. The linear structure is defined by
$\overline{\lambda x+y}=\bar \lambda\bar x+\bar y$. $\bar X$
becomes a left $A$-module when we set
$$
a\cdot\bar x=\overline{xa^*}
$$
and the $A$-valued inner product is
$$
\left<\bar x,\bar y \right>=\left< x,y \right>.
$$
This makes $\bar X$ a left Hilbert $C^*$-module over $A$.

>From how on, unless we say otherwise, all Hilbert $C^*$-modules
are assumed to be right modules and full.

A bounded module map $T:X\to X$ (where $X$ is a Hilbert $C^*$-
module) is said to be adjointable if there exists a map
$T^*:X\to X$ with $\left< Tx,y \right> =\left<x,T^*y \right>$ for all $x,y$
in $X$. The set of all adjointable maps on $X$  is a $C^*$-algebra (with respect
to the operator norm) and is denoted $\BB(X)$.

Given $X$ and $Y$ in $X$ we can define an adjointable operator
$x\otimes y^*\in\BB(X)$ by
$$
x\otimes y^*(z)=x\left<y,z \right>.
$$
(Another notation frequently used for this operator is $\theta_{x,y}$). The
$C^*$-subalgebra generated by these operators will be written $\KK(X)$.
Elements of $\KK(X)$ are sometimes referred to as ``compact operators''. If $H$
is a Hilbert space, viewed as a $C^*$-module over $\CC$, then $\KK(X)=K(H)$,
the algebra of compact operators on $H$. In general $\KK(X)\ne\BB(X)$.

Given a Hilbert $C^*$-module $X$ over $A$ one can form
$$
{\cal L}(X)=
\begin{pmatrix}
\KK(X) & X\\
\bar X & A
\end{pmatrix}.
$$
Then ${\cal L}(X)$ is a $^*$-algebra with product and involution defined by
$$
\begin{pmatrix}
K_1 & x_1\\
\bar y_1 & a_1
\end{pmatrix}
\begin{pmatrix}
K_2 & x_2\\
\bar y_2 & a_2
\end{pmatrix}
=
\begin{pmatrix}
K_1K_2+x_1\otimes y^*_2 & K_1x_2+x_1a_2\\[0.2cm]
\overline{K_2y_1}+a_1\cdot\bar y_2 & \left<y_1,x_2\right>+a_1a_2
\end{pmatrix}
$$
and
$$
\begin{pmatrix}
K & x\\
\bar y & a
\end{pmatrix}^*
=
\begin{pmatrix}
K^* & y\\
\bar x & a^*
\end{pmatrix}.
$$
There is also a natural action of ${\cal L}(X)$ on $X\oplus A$
which defines a norm on ${\cal L}(X)$ making it a $C^*$-algebra.
We shall refer to ${\cal L}(X)$ as the {\it linking algebra} of $X$.

A (right) Hilbert $C^*$-module $X$ over $A$ is said to be {\it selfdual}
if for every $A$-module map
$$
f:X\to A
$$
there is some $y\in X$ such that $f(x)=\left<y,x \right>$. Suppose now that $X$ is
a selfdual Hilbert $C^*$-module over a von Neumann algebra $M$. Then $X$ is a dual
Banach space (i.e. there is a Banach space $X_*$ such that $X=(X_*)^*$) and $\BB(X)$
is a von Neumann algebra. (See \cite[Proposition~3.8 and Proposition~3.10]{P}.)
In this case we set
$$
{\cal L}_w(X)=
\begin{pmatrix}
\BB\,(X) & X\\
\bar X & M
\end{pmatrix}.
$$
This is then a von Neumann algebra which we call the {\it linking von~Neumann
algebra} of $X$. (See~\cite{B2}).

For more about Hilbert $C^*$-modules see~\cite{L2}, \cite{RW} and \cite{P}.

%                               Page 12?????????

%                                  SECTION  2

\section{Isometries of selfdual modules.}

The main theorem in this section is the following.

%                             Theorem 2.1

\begin{theorem}
Let $M,N$ be von Neumann algebras and $p\in N, \ q\in N$ be projections such
that each of the projections $p, \ I-p, \ q$ and $I-q$ has central support equal to
$I$. Let $T:pM(I-p)\to qN(I-q)$ be a surjective linear isometry. Then there are
central projections $e_1,e_2$ in $M,f_1,f_2$ in $N$ with $e_1+e_2=I$ and $f_1+f_2=I$
and there are maps
\be
\Psi:e_1Me_1\to f_1Nf_1\nonumber\\
\Phi:e_2Me_2\to f_2Nf_2\nonumber
\ee
satisfying
\begin{description}
\item (1) $\Psi$ is a (surjective) $^*$-isomorphism\\
        \hspace*{-0.2cm}  $\Phi$ is a (surjective) $^*$-antiisomorphism.

\item(2) For $x\in pM(I-p)$,
$$
T(x)=\Psi(e_1 x e_1)+\Phi(e_2 x e_2).
$$

\item(3) $\Psi(e_1p)=f_1q$ \ and \ $\Phi(e_2p)=f_2(I-q)$.
\end{description}
\end{theorem}

The proof will be divided into several lemmas and propositions. The final
arguments can be found following Corollary~2.23.

Note first that both $pM(I-p)$ and $qN(I-q)$ have a structure of a $JB^*$-triple
with $\{x,y,z\}=\,\frac 12 (xy^*z+zy^*x)$.

Since linear isometries preserve the triple product we have the following result
which is well known (see~\cite[Proposition~5.5]{Ka} or \cite[Theorem~4]{H}).

%                                 LEMMA 2.2

\begin{lemma}
For $x,y,z\in pM(I-p)$,
$$
T(xy^*z+zy^*x)=T(x)T(y)^*T(z)+T(z)T(y)^*T(x)
$$
and, in particular,
$$
T(xy^*x)=T(x)T(y)^*T(x).
$$
\end{lemma}

Note that an element $u\in pM(I-p)$ is a partial isometry if and only if it is
a tripotent (i.e. $\{u,u,u\}=u$); hence $T(u)$ is also a partial isometry.

>From now on we assume that $M,N,p,q$ and $T$ are as in assumptions of Theorem~2.1.

%                        LEMMA 2.3

\begin{lemma}
Suppose $\{z_\alpha\}_{\alpha\in\Lambda}$ is an orthogonal family of central
projections in $M$ with $\sum z_\alpha=I$.
Then there is an orthogonal family $\{c_\alpha\}_{\alpha\in\Lambda}$ of central
projections in $N$ with $\sum c_\alpha=I$ and such that, for every
$\alpha\in\Lambda$,
$$
T\bigl(z_\alpha pM(I-p)z_\alpha\bigr)=c_\alpha qN(I-q)c_\alpha.
$$
\end{lemma}

%                      PROOF
\noindent
{\bf Proof.}~~Suppose $M\subseteq B(H)$ and $N\subseteq B(K)$ for Hilbert spaces
$H,K$.
For each $\alpha\in\Lambda$ write
\be
&&c^1_\alpha=\bigl[T\bigl(z_\alpha pM(I-p)z_\alpha\bigr)K \bigr]\nonumber\\
&&c^2_\alpha=\bigl[T\bigl(z_\alpha pM(I-p)z_\alpha\bigr)^*K \bigr]\nonumber
\ee
(where $[S]$ for a subspace $S\subseteq K$ is the projection onto the closure of $S$).
Clearly $c^1_\alpha,\ c^2_\alpha\in N$. Also $c^1_\alpha\le q, \ c^2_\alpha\le I-q$
(in particular $c^1_\alpha c^2_\alpha =0$).

If $u\in T\bigl(z_\alpha pM(I-p)z_\alpha\bigr)$ and
$v\in T\bigl(z_\beta pM(I-p)z_\beta\bigr)$ are partial isometries and
$\alpha\ne\beta$ in $\Lambda$ then $T^{-1}(u),T^{-1}(v)$ are orthogonal; hence $u$
and $v$ are orthogonal. Since $pM(I-p)$ is spanned by its partial isometries we find
that
$$
T\bigl(z_\beta pM(I-p)z_\beta\bigr)^*T \bigl(z_\alpha pM(I-p)z_\alpha\bigr)=0.
$$
It follows that
$$
c^1_\alpha c^1_\beta =0 \ \ \ \alpha\ne\beta
$$
and similarly
$$
c^2_\alpha c^2_\beta =0 \ \ \ \alpha\ne\beta.
$$
We now set $c_\alpha=c^1_\alpha + c^2_\alpha \in N$ and we find that
$\{c_\alpha\}_{\alpha\in\Lambda}$ is an orthogonal family of projections in $N$.
Now suppose $y\in qN(I-q)$ satisfies $c^1_\alpha y=0$ for all $\alpha\in\Lambda$.
Then $T(z_\alpha x)^* y=0$ for all $x\in pM(1-p)$. In particular
$$
yT\bigl(z_\alpha T^{-1}(y)\bigr)^*y=0.
$$
Applying $T^{-1}$ we get
$$
T^{-1}(y)\bigl(z_\alpha T^{-1}(y)\bigr)^*T^{-1}(y)=0.
$$
Hence $z_\alpha T^{-1}(y)\bigl(z_\alpha T^{-1}(y)\bigr)^*z_\alpha T^{-1}(y)=0$;
thus $z_\alpha T^{-1}(y)=0$ for all $\alpha\in\Lambda$. As $\sum z_\alpha=I, \
T^{-1}(y)=0$ and, consequently, $y=0$.

This proves that $\sum c^1_\alpha=q$ (recall that the central support of
$I-q$ is $I$). Similarly $\sum c^2_\alpha=I-q$. Hence $\sum c_\alpha=I$.
To show that each $c_\alpha$ is in the center of $N$ it suffices to show that,
for $\alpha\ne\beta$
$$
c_\alpha N c_\beta=0.
$$
To show that $c^2_\alpha N c^1_\beta=0$ we fix $y\in M$, and $x_2,x_1\in pM(1-p)$
and compute
\be
&&T(z_\alpha x_1)y^* T(z_\beta x_2)+T(z_\beta x_2)y^* T(z_\alpha x_1)\nonumber\\
&&=T\left(z_\alpha x_1T^{-1}((qy(I-q))\bigr)^*z_\beta x_2+
z_\beta x_2T^{-1}\bigl((qy(I-q))\bigr)^*z_\alpha x_1 \right)\nonumber\\
&& =T(0+0) \ \ ({\rm as} \ z_\alpha z_\beta=0).\nonumber
\ee
Since $c^1_\alpha T(z_\alpha x_1)=T(z_\alpha x_1)$ while $c^1_\alpha T(z_\beta x_2)=0$,
we get
$$
T(z_\alpha x_1)y^* T(z_\beta x_2)=0
$$
for all $y\in N, \ x_1,x_2\in pM(I-p)$. Hence
$$
c^2_\alpha N c^1_\beta=0.
$$
Now we turn to the proof of $c^1_\alpha Nc^1_\beta=0$. We have to show that
$c^1_\alpha qNq c^1_\beta=0$ and, since $qNq$ is the $\sigma$-weak closure of span
of $\{ab^*:a,b\in qN(I-q)\}$, we need to show: $c^1_\alpha ab^* c^1_\beta=0$ for
all such $a,b$.

Write $a=T(d)$ and $b=T(g)$ and then compute, for $x_1,x_2\in pM(I-p)$,
\be
&& T(z_\alpha x_1)^* T(d)T(g)^* T(z_\beta x_2)\nonumber\\
&&=T(z_\alpha x_1)^* \bigl[T(d)T(g)^* T(z_\beta x_2)+
T(z_\beta x_2)T(g)^*T(d)\bigr]\nonumber\\
&&=T(z_\alpha x_1)^* T\bigl(dg^*z_\beta x_2+z_\beta x_2 g^*d\bigr)\nonumber\\
&&=T(z_\alpha x_1)^* c^1_\alpha c^1_\beta T\bigl(z_\beta dg^*x_2+z_\beta x_2g^*d\bigr)=0.\nonumber
\ee
Hence $c^1_\alpha Nc^1_\beta=0$.

Thus $c_\alpha Nc_\beta=0$ for all $\alpha\ne\beta$; i.e., $c_\alpha\in Z(N)$. \ \ \
$\blacksquare$

%                              LEMMA 2.4
\begin{lemma}
Let $\{z_\alpha\}_{\alpha\in\Lambda}$ be an orthogonal family of central projections
of $M$ with $\sum z_\alpha=I$. Let $\{c_\alpha\}$ be as Lemma~2.3. Suppose that, for
every $\alpha\in\Lambda$, Theorem~2.1 holds for $z_\alpha M,c_\alpha N$ and the
restriction of $T$ to $z_\alpha pM(I-p)z_\alpha$, in place of $M,N$ and $T$.
Then it holds for $M,N,T.$
\end{lemma}
%                           PROOF
\medskip

\noindent
{\bf Proof.}~~From the conclusion of Theorem~2.1, applied to $z_\alpha M$ and
$c_\alpha N$, we get projections $e_{1,\alpha},\, e_{2,\alpha}, \, f_{1,\alpha}, \,
f_{2,\alpha}$ and maps $\Psi_\alpha,\Phi_\alpha$. Setting $e_i=\sum e_{i,\alpha}, \
f_i=\sum f_{i,\alpha}, \ \Psi=\sum\oplus\Psi_\alpha$ and $\Phi=\sum\oplus\Psi_\alpha$
we obtain the conclusion of the theorem for $M,N$. \ \ \ \ $\blacksquare$

%                              LEMMA 2.5
\begin{lemma}
There is an orthogonal family of central projections $\{z_\alpha:\alpha\in\Lambda\}$
in $M$ with $\sum z_\alpha=I$ such that, for each $\alpha\in\Lambda$, either
\begin{description}
\item(1) $z_\alpha p$ and $z_\alpha(I-p)$ are abelian projections in $M$ or
\item(2) There is a family $\{u_i:i\in I\}$, of cardinality $|I|\ge 2$, of partial
isometries in $z_\alpha pM(I-p)z_\alpha$ satisfying

\begin{description}
\item(i)~~~$u^*_iu_i=u^*_ju_j$ for all $i,j\in I$.
\item(ii)~~$u_iu^*_iu_ju^*_j=0$  for all $i\ne j$ in $I$.
\item(iii)~$\sum u_iu^*_i=pz_\alpha$
\end{description}
or
\item(3) There is a family $\{u_i:i\in I\}$, of cardinality $|I|\ge 2$, of partial
isometries in $z_\alpha pM(1-p)z_\alpha$ satisfying

\begin{description}
\item(i$'$)~~~$u_iu^*_i=u_ju^*_j$ for all $i,j\in I$.
\item(ii$'$)~~$u^*_iu_iu^*_ju_j=0$ for all $i\ne j$ in $I$.
\item(iii$'$)~$\sum u^*_iu_i=z_\alpha(I-p)$.
\end{description}
\end{description}
\end{lemma}

%                                 PROOF

\noindent {\bf Proof.}~~Since $M$ can be written as a direct sum
of algebras of different types we can deal with each type
separately. Recall that for a projection $g, \ c(g)$  is its
central support.

%                               CASE 1

\noindent
\underline{Case~1:}~~$M$ is of type III.\\
Then we can write $p=p_1+p_2$ with $p_1\sim p_2$. Since
$c(p_1)=c(p_2)=c(p)=I=c(I-p)$ and $p_1,I-p$ are both properly infinite,
$p_1\sim I-p \ \ i=1,2$. Hence there are $u_1,u_2$ in $M$ such that
$u_iu^*_i=p_i, \ u^*_iu_i=I-p, \ i=1,2$, and we are done.
\medskip

%                               CASE 2

\noindent
\underline{Case~2:}~~$M$ is of type I.\\
In this case there is an abelian projection $e_1\in M$ with $c(e_1)=I$. Since
$c(e_1)\le c(I-p) \ (=I)$, we have $e_1\precsim I-p$ (\cite[Proposition~6.4.6]{KR})
and, thus, there is an abelian projection $e\le I-p$ with $c(e)=I$. It now follows
\cite[Corollary~6.5.5]{KR} that there is  a family $\{q_j:1\le j\le\infty\}$ of
central projections with $\sum q_j=I$ and such that $q_jp$ is the sum of $j$ equivalent
abelian projections $q_jp=\sum\limits^j_{i=1} q_jp_i$. As
$c(q_jp_i)=c(q_jp)=q_j=c(q_je)$ we have $q_jp_i\sim q_je$ for all $i\le j$. Hence
for each algebra $q_jMq_j$ with $j\ge 2$ statement~(2) holds. It is left to deal
with the case $q_1$. So we  assume now that $p$ is abelian. If $I-p$ is not abelian
we can use a similar argument to the one above (interchanging the roles of $p$ and
$I-p$) and get statement~(3). We are left with the case where both $p$ and $I-p$
are abelian projections and this is (1).

%                               CASE 3

\noindent
\underline{Case~3:}~~$M$ is of type II.\\
By splitting $M$ using central projections we can assume that each of the projections
$p$ and $I-p$ are either finite or properly infinite.

If $I-p$ is properly infinite we can argue as in the type~III case:
$p=p_1+p_2, \linebreak
 p_1\sim p_2\precsim I-p$ and, thus, there are $u_1,u_2$ with
$p_i=u_iu^*_i, \ u^*_1u_1=u^*_2u_2\le I-p$; hence statement~(2) holds. If $p$ is
properly infinite we reverse the roles of $p$ and $I-p$ and get statement~(3).
So we assume that both $p$ and $I-p$ are finite (thus we are in the type~II$_1$ case).
In this case we let $\Delta$ be the center-valued dimension function, defined on the
projections of $M$ with range equal to the set of all positive operators in the unit
ball of the center (see~\cite[\S\,8.4]{KR}).

For every $j\ge 2$ we can let $q_j$ be the maximal central projection satisfying
$\frac 1j\,q_j\Delta(p)\le q_j\Delta(I-p)$, and $q_0=I-\vee q_j$. But, for every $j\ge 2$,
$$
q_0\Delta(I-p)\le q_0\,\frac 1j\,\Delta(p)\le\,\frac 1j\, q_o.
$$
Hence $q_0\Delta(I-p)=0$ and $\Delta(q_0(I-p))=0$ implying that $q_0(I-p)=0$.
But $c(q_0(I-p))=q_0c(I-p)=q_0$; hence $q_0=0$ and
$I=\sum\limits^\infty_{j=2}(q_j-q_{j-1})$ (setting $q_1=0$). Given $j\ge2$,
$$
\frac 1j\,(q_j-q_{j-1})\Delta(p)\le (q_j-q_{j-1})\Delta(I-p).
$$
Restricting our attention to the algebra $(q_j-q_{j-1})M(q_j-q_{j-1})$ we can write
$\frac 1j\,\Delta(p)\le(I-p)$. Thus we can write $p$ as a sum of $j$ equivalent
subprojections $p=\sum p_i$ with $p_i\precsim I-p$; hence $p_i\sim e\le I-p$ for all
$i\le j$. This shows that, in this case, (2) holds. \ \ \ \ $\blacksquare$

%                             LEMMA 2.6

\begin{lemma}
If $M,N,T$ are as in Theorem~2.1 and, in addition, $p$ and $I-p$ are abelian
projections in $M$, then Theorem~2.1 holds.
\end{lemma}

%                                 PROOF

\noindent {\bf Proof.}~~Since $p$ and $I-p$ are abelian projections, $M$ is
$^*$-isomorphic to $M_2\otimes{\cal A}$ where $\cal A$ is an abelian von~Neumann
algebra and $M_2$ is the algebra of $2\times 2$ complex matrices. We assume now
that $M=M_2\otimes{\cal A}$. Write $u=e_{12}\otimes I$ (where $\{e_{ij}\}$ are the
matrix units in $M_2$) and $v=T(u)$. Given $a\in{\cal A}$ we have
$$
T(e_{12}\otimes a)=T\bigl((e_{12}\otimes I)(e_{12}\otimes a^*)^*(e_{12}\otimes I)\bigr)
=vT(e_{12}\otimes a^*)^*v.
$$
Then $qN(1-q)=T(e_{12}\otimes{\cal A})=vT(e_{12}\otimes{\cal A})^*v$.
Hence $q\le vv^*$ and $I-q\le v^*v$. But $v\in qN(1-q)$ so that $q=vv^*, \ I-q=v^*v$.

Now write
$$
\psi_{1,1}(a)=T(e_{12}\otimes a)v^*\in qN(1-q)v^*=qNq, \ \ a\in{\cal A}.
$$
$\psi_{1,1}$ maps $I$ into $vv^*=q$ and it is an isometry  onto the von Neumann
algebra $qNq$. By \cite[Theorem~10]{K} this map is a $^*$-isomorphism (using the
fact that $\cal A$ is abelian). Hence $qNq$ is abelian. Similarly one sees that $I-q$
is also an abelian projection.

We now have, for $a,b,c\in{\cal A}$,
\be
&&T(e_{12}\otimes a)v^*T(e_{12}\otimes b)=T(e_{12}\otimes a)v^*T(e_{12}\otimes b)v^*v\nonumber\\
&&=\psi_{11}(a)\psi_{11}(b)v=\psi_{11}(ab)v=T(e_{12}\otimes ab)\nonumber
\ee
\be
&&T(e_{12}\otimes b)T(e_{12}\otimes c^*)^*=\psi_{11}(b)v[\psi_{11}(c^*)v]^*=\psi_{11}(b)vv^*\psi_{11}(c^*)^*\nonumber\\
&&=\psi_{11}(b)\psi_{11}(c)=T(e_{12}\otimes bc)v^*\nonumber
\ee
and, using similar identities we see that the map
$$
\psi:M\to N
$$
defined by
$$
\psi
\begin{pmatrix}
a & b\\
c & d
\end{pmatrix}
=
\begin{pmatrix}
T(e_{12}\otimes a)v^*  & T(e_{12}\otimes b)\\[0.2cm]
T(e_{12}\otimes c^*)^* & v^*T(e_{12}\otimes d)
\end{pmatrix}
$$
is a $^*$-isomorphism of $M$ onto $N$ extending $T$. This completes the proof
of Theorem~2.1 in this case. (Here $e_1=I, \ e_2=0$.) \ \ \ \ $\blacksquare$

>From now on, in this section, we assume that condition~(2) of Lemma~2.5 holds (with
$z_\alpha=I$).

We fix the family $\{u_i\}$ as in Lemma~2.5 and write $v_i$ for $T(u_i)\in qN(1-q)$.
Then $v_i$ is a partial isometry and we write
$$
r_i=v_iv^*_i(\le q), \ d_i=v^*_iv_i(\le I-q).
$$
Now fix $i\ne j$. We wish to study the relative position of $v_i$ and $v_j$.

We have
\be
&&v_i=T(u_iu^*_iu_i)=T(u_iu^*_ju_j)=T(u_iu^*_ju_j+u_ju^*_ju_i)\nonumber\\
&&=v_iv^*_jv_j+v_jv^*_jv_i=v_id_j+r_jv_i.\nonumber
\ee
It then follows that $r_jv_id_j=0$.
Since we can interchange $i$ and $j$ we get
\be
v_i=v_id_j+r_jv_i \ \ \ {\rm and} \ \ \ r_jv_id_j=0\nonumber\\
v_j=v_jd_i+r_iv_j \ \ \ {\rm and} \ \ \ r_iv_jd_i=0\nonumber
\ee
and, thus,
$$
r_i=v_iv^*_i=(v_id_j+r_jv_i)(d_jv^*_i+v^*_ir_j)=
v_id_jv^*_i+r_jr_ir_j.
$$
But then $r_jr_ir_j\le r_i$ and, consequently $(I-r_i)r_jr_ir_j(I-r_i)=0$ which
implies that $r_ir_j=r_ir_jr_i$ and $r_ir_j=r_jr_i$. Similar analysis works for
$d_i,d_j$ and we find that
$$
r_ir_j=r_jr_i, \ d_id_j=d_jd_i.
$$
The computation above shows now that
$$
r_i=v_id_jv^*_i+r_ir_j
$$
and similarly
$$
d_i=v^*_ir_jv_i+d_id_j.
$$
We summarize as follows.

%                              LEMMA 2.7
\begin{lemma}
With the notation and assumptions above, for $i\ne j$,
\begin{description}
\item(1)~~~$v_i=v_id_j+r_jv_i=(I-r_j)v_id_j+r_jv_i(I-d_j)$
\item(2)~~~$v_j=v_jd_i+r_iv_j=(I-r_i)v_jd_i+r_iv_j(I-d_i)$
\item(3)~~~$r_jv_id_j=r_iv_jd_i=0$
\item(4)~~~$r_ir_j=r_jr_i, \ d_id_j=d_jd_i$
\item(5)~~~$v_id_jv^*_i=r_i(I-r_j)$\\
\end{description}
\vspace{-0.5cm}

\noindent
and
\begin{description}
\item(6)~~~$v^*_ir_jv_i=d_i(I-d_j)$.\ \ \ \ $\blacksquare$
\end{description}
\end{lemma}

%                              LEMMA 2.8
\begin{lemma}
With the notation and assumption above we have for every triple $\{i,j,k\}$ of
different indices,
\begin{description}
\item(1)~~~$d_id_j(I-d_k)=0$
\end{description}

\noindent
and
\begin{description}
\item(2)~~~$r_ir_j(I-r_k)=0.$.
\end{description}
Consequently, if we write $r$ for $\vee\{r_i:i\in I\}$ and $r_0$ for
$\wedge\{r_i:i\in I\}$, then $\{r_i-r_0:i\in I\}$ is an orthogonal family
of projections with sum equal to $r-r_0$. Similar statement holds for $\{d_i-d_0\}$.

\end{lemma}

%                                 PROOF

\noindent {\bf Proof.}~~For every $a\in u_iMu^*_i$ we have
$$
au_iu^*_ju_k+u_ku^*_jau_i=0 \ \ \ ({\rm as} \ \ u^*_ju_k=0=u^*_ju_i).
$$
Thus
$$
T(au_i)v^*_jv_k+v_kv^*_j T(au_i)=0\eqno{(*)}
$$
Now set
$$
a=u_iu^*_iT^{-1}\bigl(r_j(I-r_k)v_i\bigr)u^*_i\in u_iMu^*_i.
$$
Then
\be
T(au_i)&=&T\bigl(u_iu^*_iT^{-1}(r_j(I-r_k)v_i)u^*_iu_i\bigr)\nonumber\\
       &=&v_iT\bigl(u_iT^{-1}(r_j(I-r_k)v_i)^* u_i\bigr)^*v_i\nonumber\\
       &=&v_i\bigl[v_i(v^*_ir_j(I-r_k))v_i\bigr]^*v_i\nonumber\\
       &=&r_ir_j(I-r_k)v_i.\nonumber
\ee
Using $(*)$ we have
$$
r_ir_j(I-r_k)v_iv^*_jv_k+v_kv^*_jr_ir_j(I-r_k)v_i=0.
$$
Multiplying on the left by $v_jv^*_k$ we get
$$
v_jv^*_kr_ir_j(I-r_k)v_iv^*_jv_k+v_jd_kv^*_jr_ir_j(I-r_k)v_i=0.
$$
As $v^*_k(I-r_k)=0$, the first term vanishes. From Lemma~2.7 we know that
$v_jd_kv^*_j=r_j(I-r_k)$; hence $r_j(I-r_k)r_ir_j(1-r_k)v_i=0$.

It follows that $r_ir_j(I-r_k)=0$. Statement~(1) is proved similarly and the final
statement of the lemma follows immediately. \ \ \ \ $\blacksquare$

%                              LEMMA 2.9
\begin{lemma}
\begin{description}

\item(1) For $a\in pM(I-p)$ and a partial isometry $u\in pM(I-p)$,
$$
T(uu^*au^*u)=vv^*T(a)v^*v
$$
where $v=T(u)$.

\item(2) For a partial isometry $u\in pM(I-p)$ with $T(u)=v$,
$$
T(u^*Mu^*u)=vv^*Nv^*v.
$$

\item(3) If $x,y\in pM(I-p)$ and $x^*y=yx^*=0$ then
$$
T(x)^*T(y)=0=T(y)T(x)^*.
$$
\end{description}
\end{lemma}

%                                 PROOF

\noindent {\bf Proof.}\\
\begin{description}
\vspace{-0.7cm}

\item(1) $T(uu^*au^*u)=vT(ua^*u)^*v=v[vT(a)^*v]^*v=vv^*T(a)v^*v.$

\item(2) From~(1) it follows that
$$
T(uu^*Mu^*u)\subseteq vv^*Nv^*v.
$$
Applying the same argument to $T^{-1}$ we get equality.

\item(3) Let $x,y\in pM(I-p)$ satisfy $x^*y=yx^*=0$. Let $x=u_1|x|$ be the polar
decomposition of $x$ and $y=u_2|y|$ be the one for $y$. Then $u^*_1u_2=u_2u^*_1=0$.\linebreak
If we write $v_i=T(u_i)$ then
$$
0=T(u_1u^*_1u_2+u_2u^*_1u_1)=v_1v^*_1v_2+v_2v^*_1v_1
$$
and
$$
0=T(u_1u^*_2u_1)=v_1v^*_2v_1.
$$
hence
\be
v^*_1v_2&=&v^*_1v_1v^*_1v_2+(v^*_1v_2v^*_1)v_1\nonumber\\
        &=&v^*_1(v_1v^*_1v_2+v_2v^*_1v_1)=0.\nonumber
\ee
Similarly $v_2v^*_1=0$.

Since  $x\in u_1u^*_1Mu^*_1u_1, \ y\in u_2u^*_2Mu^*_2u_2$, part~(1) shows that
$$
T(x)\in v_1v^*_1Nv^*_1v_1, \ T(y)\in v_2v^*_2Nv^*_2v_2,
$$
As $v_2v^*_1=v^*_1v_2=0$,
$$
T(x)^*T(y)=0=T(y)T(x)^*. \ \ \ \  \ \ \blacksquare
$$
\end{description}

%                              LEMMA 2.10
\begin{lemma}
If $u$ is a partial isometry in $pM(I-p), \ v=T(u)$ and $e$ is a projection
satisfying $e\le vv^*$ then there is a projection $e_0\le uu^*$ with
$T^{-1}(ev)=e_0u$.
\end{lemma}

%                                 PROOF

\noindent {\bf Proof.}~~Write $v'=ev,v''=(vv^*-e)v$. Both are partial isometries
and they satisfy $v'v''{^*}=v''v'{^*}=0$. Hence (Lemma~2.9), the partial isometries
$u'=T^{-1}(v')$ and $u''=T^{-1}(v'')$ satisfy $0=u'u''{^*}=u''u'{^*}$. Since
$u'+u''=T^{-1}(v'+v'')=T^{-1}(v)=u$ the conclusion follows. \ \ \ \ $\blacksquare$

%                              LEMMA 2.11
\begin{lemma}
For all $i,j\in I$,  \ $r_j$ commutes with the elements in $r_iNr_i$.
\end{lemma}

%                                 PROOF

\noindent {\bf Proof.}~~For $i=j$ it is clear, so assume $i\ne j$.
Since $r_iNd_i=T(u_iu^*_iMu^*_iu_i)$  (Lemma~2.9\,(2)),
$$
r_iNr_i=r_iNd_iv^*_i=T(u_iu^*_iMu^*_iu_i)v^*_i.
$$
So fix $x=u_iu^*_ixu^*_iu_i$ and compute (using the fact that $u_ju^*_ju_iu^*_i=0$ and
$u^*_ju_j=u^*_iu_i$).
\be
T(x)&=&T(xu^*_iu_i)=T(xu^*_ju_j+u_ju^*_jx)=T(x)v^*_jv_j+v_jv^*_jT(x)\nonumber\\
    &=&T(x)d_j+r_jT(x).\nonumber
\ee
Hence
$$
T(x)=(I-r_j)T(x)d_j+r_jT(x)(I-d_j).
$$
\be
T(x)v^*_ir_j&=&T(x)v^*_ir_jv_iv^*_i=T(x)d_i(I-d_j)v^*_i\nonumber\\
            &=&r_jT(x)d_i(1-d_j)v^*_i.\nonumber
\ee
Hence, for every $y\in r_iNr_i, \ yr_j=r_jyr_j$ and the claim follows. \ \ \ \ $\blacksquare$

Our next objective is to show that, for $i\ne j$ and for $x,y,z$ in
$u_iMu_i, \linebreak
 r_i(I-r_j)T(xy^*z)=r_i(I-r_j)T(x)T(y)^*T(z)$.  This will be proved in
Proposition~2.13.

We first consider the map
$$
\vr:u_iu^*_iMu_iu^*_i\to v_iv^*_iNv_iv^*_i
$$
defined by
$$
\vr(u_iu^*_ixu_iu^*_i)=T(u_iu^*_ixu_i)v^*_i.
$$
Since $T(u_iu^*_iMu_i)=r_iNd_i$ (Lemma~2.9), the map $\vr$ is a surjective isometry
from the von Neumann algebra $u_iu^*_iMu_iu^*_i$ onto the von Neumann algebra
$r_iNr_i$ that is unital $\bigl(\vr(u_iu^*_i)=r_i\bigr)$.
By~\cite[Theorem~10]{K} there are central projections $g,h$ in $r_iNr_i$ and central
projections $g_0,h_0$ in $u_iu^*_iMu_iu^*_i$ with $g+h=r_i, \ g_0+h_0=u_iu^*_i$ and
such that $\vr$, restricted to $g_0Mg_0$, is a $^*$-isomorphism onto $gNg$ and
$\vr$, restricted to $h_0Mh_0$, is a $^*$-antiisomorphism of $h_0Mh_0$ onto $hNh$.

%                              LEMMA 2.12
\begin{lemma}
With the notation above, $h(I-r_j)$ is an abelian projection in $r_iNr_i$.
\end{lemma}

%                                 PROOF

\noindent
{\bf Proof.}~~Since $h\in r_iNr_i$ it follows from Lemma~2.11 that
$h(I-r_j)$ is a projection in $r_iNr_i$. Write $c=h(I-r_j)$. To show that $cNc$
is abelian it suffices to show that one cannot find in $cNc$ projections $e_1,e_2$
that are equivalent (in $cNc$) and orthogonal (i.e. $e_1e_2=0$). Assume, by negation,
that there are such projections. Then there is a partial isometry $w\in cNc$ with
$$
ww^*=e_2 \ \ \ w^*w=e_1.
$$
Write
$$
t_1=e_1v_i, \ t_2=v_jv^*_ie_1v_i, \ t_3=wv_i
$$
and set $s_i=T^{-1}(t_i)$. Then $t_i$ and $s_i$ are partial isometries. We have
$t^*_1t_1=v^*_ie_1v_1, \ t^*_2t_2=v^*_ie_1v_iv^*_jv_jv^*_ie_1v_i$ and
$t^*_3t_3=v^*_iw^*wv_i=v^*_ie_1v_i$.
Since
$$
v^*_ie_1v_i\le v^*_icv_i\le v^*_i(I-r_j)v_i=d_id_j\le d_j,
$$
(using Lemma~2.7), we have
$$
t^*_it_i=v^*_ie_1v_i \ \ \ \ i=1,2,3.
$$
Also, $t_1t^*_1=e_1, \ t_2t^*_2=v_jv^*_ie_1v_iv^*_j\le r_j$ and
$t_3t^*_3=wv_iv^*_iw^*=e_2\le I-r_j$. Hence $\{t_it^*_i:i=1,2,3\}$ is an orthogonal
set.

By Lemma~2.10, there are projections $c_1\le u_iu^*_i$ and $c_2\le u_ju^*_j$ with
$$
s_1=c_1u_i, \ s_2=c_2u_j.
$$
Now, by the definition of $h$, the map
$$
\vr(x)=T(xu_i)v^*_i, \ x\in h_0Mh_0
$$
is a $^*$-antiisomorphism. Hence
\be
T(s_3s^*_1s_1)&=&\vr(s_3s^*_1s_1u^*_i)v_i=
\vr\bigl((s_3u^*_i)(u_is^*_1)(s_1u^*_i)\bigr)v_i\nonumber\\
              &=&\vr(s_1u^*_i)\vr(s_1u^*_i)^*\vr(s_3u^*_i)v_i\nonumber\\
              &=&\vr(s_1u^*_i)v_iv^*_i\vr(s_1u^*_i)^*\vr(s_3u^*_i)v_i\nonumber\\
              &=&T(s_1)T(s_1)^*T(s_3)=t_1t^*_1t_3=0.\nonumber
\ee
Hence $s_3s^*_1=0$.

We can use Lemma~2.8 and Lemma~2.7 for $\{t_1,t_2,t_3\}$ and $T^{-1}$ in place of
$\{v_i,v_j,v_k\}$ and $T$ (since they also have pairwise orthogonal ranges and the
same initial space). Since $s^*_3s_3s^*_1s_1=0$ (by the computation above) it follows
from Lemma~2.8\,(1) that
$s^*_2s_2s^*_1s_1=s^*_2s_2s^*_3s_3s^*_1s_1=0$ and, similarly, $s^*_2s_2s^*_3s_3=0$.
So that $\{s^*_is_i\}$ is an orthogonal family. By Lemma~2.7\,(5) (applied to the
situation here) we get for $i\ne j$ in $\{1,2,3\}$,
$s_is^*_i(I-s_js^*_j)=s_is^*_js_js^*_i=0$.

We conclude that $s_is^*_i=s_js^*_j$ for all $i,j$ in $\{1,2,3\}$. But
$s_1s^*_1\le u_iu^*_i$ (as $s_1=c_1u_i$) and $s_2s^*_2\le u_ju^*_j\le I-u_iu^*_i$.

This is a contradiction and it completes the proof. \ \ \ \ $\blacksquare$

%                             PROPOSITION 2.13
\begin{proposition}
For $x,y,z$   in $u_iMu_i$ and $j\ne i$,
$$
r_i(I-r_j)T(xy^*z)=r_i(I-r_j)T(x)T(y)^*T(z).
$$
\end{proposition}

%                                 PROOF

\noindent
{\bf Proof.}~~Fix $x,y,z$ in $u_iMu_i$ and write
$$
r_i(I-r_j)T(xy^*z)=r_i(I-r_j)hT(xy^*z)+r_i(I-r_j)gT(xy^*z)
$$
where $g,h$ were defined above. With $\vr$ as above we have
$$
(I-r_j)hT(xy^*z)=(I-r_j)h\vr(xy^*zu^*_i)v_i.
$$
Since $\vr(xy^*zu^*_i)$ lies in $r_iNr_i$, we can use Lemma~2.11 and Lemma~2.7 \linebreak
to get $(I-r_j)\vr(xy^*zu^*_i)v_i=\vr(xy^*zu^*_i)(I-r_j)v_i=
\vr(xy^*zu^*_i)v_iv^*_i(I-r_j)v_i$.
Hence, by the definition of $h$,
$(I-r_j)hT(xy^*z)=(I-r_j)h\vr(xy^*zu^*_i)hr_i(I-r_j)v_i=
(I-r_j)h\vr(zu^*_i)\vr(yu^*_i)^*\vr(xu^*_i)hr_i(I-r_j)v_i$.
But since $(I-r_j)h$ is an abelian projection in $r_iNr_i$ (Lemma~2.12) we have
\be
(I-r_j)hT(xy^*z)&=&(I-r_j)h\vr(xu^*_i)\vr(yu^*_i)^*\vr(zu^*_i)hr_i(I-r_j)v_i\nonumber\\
                &=&(I-r_j)hT(x)v^*_iv_iT(y)^*T(z)v^*_ihr_i(I-r_j)v_i\nonumber\\
                &=&(I-r_j)hT(x)T(y)^*T(z).\nonumber
\ee
Also, from the definition of $g$,
\be
(I-r_j)gT(xy^*z)&=&(I-r_j)g\vr(xy^*zu^*_i)v_i\nonumber\\
                &=&(I-r_j)g\vr(xu^*_i)\vr(yu^*_i)^*\vr(zu^*_i)v_i\nonumber\\
                &=&(I-r_j)gT(x)T(y)^*T(z).\nonumber
\ee
This completes the proof. \ \ \ \ $\blacksquare$

We now turn to define a map
$$
\theta:pMp\to N.
$$
For it we note first that $p=\sum u_iu^*_i$ and every $x\in pMp, \
x=\sum\limits_{i,j} u_iu^*_ixu_ju^*_j$ \linebreak
($\sigma$-weakly). For every $i,j\in I$ we set
$$
\theta_{ij}:u_iu^*_iMu_ju^*_j \longrightarrow N
$$
by
$$
\theta_{ij}(u_iu^*_ixu_ju^*_j)=(r_i-r_0)T(u_iu^*_ixu_j)v^*_j
+v^*_jT(u_iu^*_ixu_j)(d_i-d_0).
$$
where $r_0=\wedge r_i, \ d_0=\wedge d_i$.

To study the map $\theta$ defined by $\{\theta_{ij}\}$ we first write
$$
\alpha_{ij}(u_iu^*_ixu_ju^*_j)=(r_i-r_0)T(u_iu^*_ixu_j)v^*_j
$$
and
$$
\beta_{ij}(u_iu^*_ixu_ju^*_j)=v^*_jT(u_iu^*_ixu_j)(d_i-d_0).
$$
Also write, for every finite subset $F\subseteq I, \ p_F=\sum_{i\in F} u_iu^*_i, \
r_F=\sum_{i\in F}r_i, \linebreak
d_F=\sum_{i\in F} d_i$ and
$$
x(F)=p_Fxp_F \ \ ({\rm for} \ x\in pMp)
$$
and $\alpha_F:p_FMp_F \to N$ is defined by
$$
\alpha_F(p_Fxp_F)=\sum_{i,j\in F}\alpha_{ij}(u_iu^*_ixu_ju^*_j).
$$
(Similarly, $\beta_F$ can be defined).
We have $\alpha_F(p_FMp_F)\subseteq (r_F-r_0)N(r_F-r_0)$ and if
$F_1\subseteq F_2, \ x\in N$,
$$
\alpha_{F_1}(p_{F_1}xp_{F_1})=(r_{F_1}-r_0)\alpha_{F_2}(p_{F_2}xp_{F_2})(r_{F_1}-r_0).
$$

%                              LEMMA 2.14
\begin{lemma}
Given a finite subset $F\subseteq I, \ \alpha_F$ is a $^*$-homomorphism of $p_FMp_F$ \linebreak
onto $(r_F-r_0)N(r_F-r_0)$ and $\beta_F$ is a $^*$-antihomomorphism  of $p_FMp_F$ onto \linebreak
$(d_F-d_0)N(d_F-d_0)$.
\end{lemma}

%                                 PROOF

\noindent
{\bf Proof.}~~We prove the properties of $\alpha_F$. The proof for $\beta_F$ is similar.

For $i,j\in I$, \ $u_iu^*_iMu_j=u_iu^*_iMu_ju^*_iu_i\subseteq u_iu^*_iMu^*_iu_i=
u_iu^*_iMu^*_ju_j\subseteq u_iu^*_iMu_j$. Hence
$u_iu^*_iMu_j=u_iu^*_iMu^*_iu_i$ and $T(u_iu^*_iMu_j)=r_iNd_i$. Hence
$T(u_iu^*_iMu_j)v^*_j=r_iNd_iv^*_j$. For $i=j$ this is $r_iNr_i$ and, thus,
$\alpha_{ii}(u_iu^*_iMu_iu^*_i)=(r_i-r_0)Nr_i=(r_i-r_0)N(r_i-r_0)$ (as $I-r_0$
commutes with $r_iNr_i$ by Lemma~2.11).

For $i\ne j$
\be
\alpha_{ij}(u_iu^*_iMu_ju^*_j)=(r_i-r_0)Nd_iv^*_j=(r_i-r_0)Nv^*_j(v_jd_iv^*_j)
=(r_i-r_0)Nv^*_j(I-r_i).\nonumber
\ee
As $r_j(I-r_i)=r_j-r_0$ (Lemma~2.8) we see that $\alpha_{ij}$ is onto
$(r_i-r_0)N(r_j-r_0)$.
This shows that $\alpha_F$ maps $p_FMp_F$ onto $(r_F-r_0)N(r_F-r_0)$.

We now show that $\alpha_F$ is a $^*$-map. For that, fix $x=u_iu^*_ixu_ju^*_j$ and
consider
$$
\alpha_{ij}(x)=(r_i-r_0)T(xu_j)v^*_j=(r_i-r_0)T(xu_j)v^*_j(r_j-r_0).
$$
If $i=j$ we have
\be
\alpha_{ii}(x)&=&(r_i-r_0)T(u_iu^*_ixu_i)v^*_i(r_i-r_0)\nonumber\\
              &=&(r_i-r_0)v_iT(x^*u_i)^*v_iv^*_i(r_i-r_0)\nonumber\\
              &=&\bigl[(r_i-r_0)T(x^*u_i)v^*_i(r_i-r_0) \bigr]^*\nonumber\\
              &=&\alpha_{ii}(x^*)^*.\nonumber
\ee
If $i\ne j$ we have $T(u_iu^*_ixu_j)=T(u_iu^*_ixu_j+u_ju^*_ixu_i)=
v_iT(x^*u_i)^*v_j+v_jT(x^*u_i)^*v_i$. Hence
$$
\alpha_{ij}(x)=(r_i-r_0)\bigl[v_iT(x^*u_i)^*v_j+v_jT(x^*u_i)^*v_i\bigr]v^*_j(r_j-r_0).
$$
Since $(r_i-r_0)v_j=0$ we have
$$
\alpha_{ij}(x)=(r_i-r_0)v_iT(x^*u_i)^*(r_j-r_0)=\alpha_{ji}(x^*)^*.
$$
This shows that $\alpha_F$ is a $^*$-map.

Finally, we shall show that $\alpha_F$ is a homomorphism. For that we fix
$$
x=u_iu^*_ixu_ku^*_k, \ y=u_ju^*_jyu_mu^*_m.
$$
If $k\ne j$ then $xy=0$. In this case $r_k(r_j-r_0)=0$ and
$\alpha_F(x)\alpha_F(y)=(r_i-r_0)T(xu_k)v^*_k(r_j-r_0)T(yu_m)v^*_m=0$.
So we suppose now that $k=j$ and then $xy=u_iu^*_ixu_ju^*_jyu_mu^*_m$. Hence
$\alpha_F(xy)=(r_i-r_0)T(xyu_m)v^*_m$.
If $i\ne j$ this is equal to
\be
&&(r_i-r_0)T\bigl((xu_j)u^*_j(yu_m)\bigr)v^*_m
=(r_i-r_0)T\bigl[(xu_j)u^*_j(yu_m)+yu_mu^*_jxu_j\bigr]v^*_m\nonumber\\
&&=(r_i-r_0)\bigl[T(xu_j)v^*_jT(yu_m)+T(yu_m)v^*_j T(xu_j)\bigr]v^*_m.\nonumber
\ee
Now $T(yu_m)=T(yu_mu^*_ju_j)\in r_jNd_j$ and $(r_i-r_0)T(yu_m)=0$. Hence
$$
\alpha_F(xy)=(r_i-r_0)T(xy_j)v^*_jT(yu_m)v^*_m=\alpha_F(x)\alpha_F(y).
$$
Now consider the case $i=j$. Then
$$
\alpha_F(xy)=(r_i-r_0)T\bigl((xu_i)u^*_i(yu_m)\bigr)v^*_m.
$$
Since $yu_m=u_iu^*_iyu_mu^*_iu_i\in u_iMu_i$ and also $xu_i,u_i$ lie in $u_iMu_i$,
we can apply Proposition~2.13 and get
$$
\alpha_F(xy)=(r_i-r_0)T(xu_i)v^*_iT(yu_m)v^*_m=\alpha_F(x)\alpha_F(y). \ \ \  \ \blacksquare
$$
It follows from Lemma~2.14 that for each finite subset $F\subseteq I$ and each $x\in pMp$
$$
\|\alpha_F(x)\|\le\|x\|.
$$
Hence, for a fixed $x\in pMp$ the net $\{\alpha_F(x):F\subseteq I\}$ is bounded and we can
find a $\sigma$-weakly convergent subset $\alpha_{F'}(x)\longrightarrow\alpha(x)$.

For every finite subset $F\subseteq I$ there is some $F'_0$ in the subnet with
$F'_0\supseteq F$. But then, for every
$F'\supseteq F'_0, \ r_F\alpha_{F'}(x)r_F=\alpha_F(x)$; hence
$$
r_F\alpha(x)r_F=\alpha_F(x)
$$
and, consequently.
$$
\alpha(x)=\sigma{\rm -weak} \lim_F \alpha_F(x).
$$
We can now conclude from Lemma~2.14 the following. (The statement for $\beta$ is
proved similarly.)

%                              COROLLARY 2.15
\begin{corollary}
There is a surjective \ $^*$-homomorphism \ $\alpha:pMp\to (r-r_0) N(r-r_0)$ and a
surjective $^*$-antihomomorphism $\beta:pMp\to (d-d_0)N(d-d_0)$ such that, for all
$i,j$,
\be
&&r_i\alpha(x)r_j=\alpha_{ij}(x)\nonumber\\
&&d_j\beta(x)d_i=\beta_{ij}(x).\nonumber
\ee
\end{corollary}

%                            LEMMA 2.16
\begin{lemma}
For $x\in pM(I-p)$ and $i,j\in I$ we have $r_i(I-r_j)T(u_ju^*_jx)=0 $.
\end{lemma}

%                              PROOF
\noindent
{\bf Proof.}~~We can assume $i\ne j$ and $x=u_ju^*_jx$. Write $x=x_1+x_2$ where
$x_1=xu^*_ju_j$ and $x_2=x(I-u^*_ju_j)$. Then $x_1\in u_jAu_j$ and, thus,
$T(x_1)\in v_jNv_j=r_jNd_j$.
Hence $(I-r_j)T(u_ju^*_jx_1)=0$. For $x_2$ we have $u_iu^*_ix_2=x_2u^*_iu_i=0$ (as
$u^*_iu_i=u^*_ju_j)$. Hence
$0=T(u_iu^*_ix_2+x_2u^*_iu_i)=v_iv^*_iT(x_2)+T(x_2)v^*_iv_i=r_iT(x_2)+T(x_2)d_i$.
We also have $u_ix^*_2u_i=0$; hence $v_iT(x_2)^*v_i=0$. Therefore
\be
r_iT(x_2)&=&r_iT(x_2)d_i+r_iT(x_2)(I-d_i)\nonumber\\
         &=&r_iT(x_2)d_i+\bigl[r_iT(x_2)+T(x_2)d_i\bigr](I-d_i)=0.\nonumber
\ee
Note that a similar argument shows that
$$
T(x_2)d_i=0. \ \ \ \ \ \blacksquare
$$

%                            LEMMA 2.17
\begin{lemma}
For $x\in u_ju^*_jM(I-p), \ y\in u_iu^*_iMu_i$,
$$
r_i(I-r_0)T(y)v^*_jT(x)=r_i(I-r_0)T(yu^*_jx).
$$
\end{lemma}

%                              PROOF
\noindent
{\bf Proof.}~~Assume first that $i\ne j$. Then it follows from Lemma~2.16 that \linebreak
$r_i(I-r_0)T(x)=0$. Also we have $xu^*_jy=0$. Hence
\be
&&r_i(I-r_0)T(yu^*_jx)=r_i(I-r_j)T(yu^*_jx+xu^*_jy)\nonumber\\
&&=r_i(I-r_j)\bigl[T(y)v^*_jT(x)+T(x)v^*_jT(y)\bigr]\nonumber\\
&&=r_i(I-r_j)T(y)v^*_jT(x).\nonumber
\ee
Now consider the case $i=j$. if $x\in u_iu^*_iMu^*_iu_i$ the result follows from
Proposition~2.13. So assume now that $x=u_iu^*_ix(I-u^*_iu_i)$. Then $xu^*_iy=0$ and
$r_iT(x)v^*_i=r_iT(x)d_iv^*_i=T(u_iu^*_ixu^*_iu_i)v^*_i=0$ (where we used Lemma~2.9).
Hence
\be
r_iT(yu^*_ix)&=&r_iT(yu^*_ix+xu^*_iy)\nonumber\\
             &=&r_i\bigl[T(y)v^*_iT(x)+T(x)v^*_iT(y)\bigr]\nonumber\\
             &=&r_iT(y)v^*_iT(x).\nonumber \ \ \ \ \ \blacksquare
\ee

%                              COROLLARY 2.18
\begin{corollary}
For every $i,j,k$ if $a=u_iu^*_iau_ju^*_j$ and $x=u_ku^*_kx(I-p)$ then
$$
(r-r_0)T(ax)=\alpha_{ij}(a)T(x).
$$
\end{corollary}

%                              PROOF
\noindent
{\bf Proof.}~~Assume first that $j\ne k$. Then $ax=0$. Also
$\alpha_{ij}(a)\in(r_i-r_0)N(r_j-r_0)$ (see the proof of Lemma~2.14) and
$(r_j-r_0)T(x)=r_j(I-r_k)T(x)=0$ by Lemma~2.16. Hence $\alpha_{ij}(a)T(x)=0$.
We now consider the case $j=k$. In this case
$\alpha_{ij}(a)T(x)=(r_i-r_0)T(au_j)v^*_jT(x)$ and Lemma~2.17 (with
$y=au_j\in u_iu^*_iMu^*_ju_j=u_iu^*_iMu^*_iu_i$) applies to give
$$
\alpha_{ij}(a)T(x)=r_i(I-r_0)T(au_ju^*_jx)=r_i(I-r_0)T(ax).
$$
As $(r-r_i)T(ax)=0$ (Lemma~2.16), we are done.\ \ \ \ $\blacksquare$

Before we conclude from the last corollary that $T$ is a module map we need the
following.

%                            LEMMA 2.19
\begin{lemma}
\
\begin{description}

\item(1) $\theta(=\alpha+\beta)$ is an injective map.

\item(2) $\alpha$ and $\beta$ are $\sigma$-weakly continuous maps on $pMp$.

\item(3) There are projections $g_1,g_2$ in $Z(pMp)$ such that
\begin{description}
\item(i)~~~$g_1+g_2=p$.
\item(ii)~~$\ker\alpha=g_2Mg_2$ and $\ker\beta=g_1Mg_1$.
\item(iii)~if $c(g_i)$ is the central support of $g_i$ in $M$, then $c(g_1)+c(g_2)=I$.
\end{description}
\end{description}
\end{lemma}

%                              PROOF
\noindent
{\bf Proof.}\\
(1)\\
Recall that for every $i,j\in I$,
\be
\theta(u_iu^*_ixu_ju^*_j)&=&(r_i-r_0)T(u_iu^*_ixu_j)v^*_j(r_j-r_0)\nonumber\\
                         &+&(d_j-d_0)v^*_jT(u_iu^*_ixu_j)(d_i-d_0).\nonumber
\ee
Since $\{r_i-r_0\}$ and $\{d_i-d_0\}$ are orthogonal families, it will suffice to
show the injectivity of $\theta_{ij} \ \bigl(=\theta|u_iu^*_iMu_ju^*_j\bigr)$
for all $i,j\in I$.

So fix $i,j\in I$ and $x=u_iu^*_ixu_ju^*_j$ such that
$$
(r_i-r_0)T(u_iu^*_ixu_j)v^*_j(r_j-r_0)=0=(d_j-d_0)v^*_jT(u_iu^*_ixu_j)(d_i-d_0).
$$
>From Lemma~2.7 we get $v_j(d_j-d_0)v^*_j=r_0$ and $v^*_j(r_j-r_0)v_j=d_0$ and we have
$$
(r_i-r_0)T(u_iu^*_ixu_j)d_0=0=r_0T(u_iu^*_ixu_j)(d_i-d_0).
$$
Also, fix $k\ne i$, and compute
$$
r_0T(u_iu^*_ixu_j)d_0=r_kT(u_iu^*_ixu_j)d_k=T(u_ku^*_ku_iu^*_ixu_ju^*_ku_k)=0
$$
(using Lemma~2.9).

Now we will show that $(r_i-r_0)T(u_iu^*_ixu_j)(d_i-d_0)=0$. It will then follow that
$T(u_iu^*_ixu_j)=0$ (as it lies in $r_iNd_i$ by Lemma~2.9) and consequently,
$u_iu^*_ixu_j=0$ implying $x=0$.

Fix $k\ne i$ and note that
$$
(r_i-r_0)T(u_iu^*_ixu_j)(d_i-d_0)=(r_i-r_0)T(u_iu^*_ixu_j)v_ir_kv^*_i
\in(r_i-r_0)Nr_ir_kv^*_i.
$$
But the last set is $\{0\}$ since $r_k$ commutes with $r_iNr_i$ (Lemma~2.11) and \linebreak
$r_k(r_i-r_0)=0$.
\medskip

\noindent
(2)\\
$\beta$ is a map onto $(d-d_0)N(d-d_0)$ and $\alpha$ is onto $(r-r_0)N(r-r_0)$.
Since $d-d_0$ and $r-r_0$ are orthogonal projections, we can view $N$ as acting on a
Hilbert space $H$ with two orthogonal subspaces $H_1=(d-d_0)(H)$ and $H_2=(r-r_0)(H)$.
Let $\tau:B(H_1)\to B(H_1)$ be a transpose map, then $\alpha\oplus\tau\circ\beta$
is a $^*$-isomorphism of $pMp$ into $B(H)$. Thus it is $\sigma$-weakly continuous and
so are its compressions to $H_1$ and $H_2$; i.e. $\alpha$ and $\beta$ are
$\sigma$-weakly continuous.

\medskip

\noindent
(3)\\
Since $\alpha,\beta$ are $\sigma$-weakly continuous their kernels are $\sigma$-weakly
closed ideals in $pMp$ and the existence of projections $g_1,g_2$ in the center of
$pMp$ and satisfying (ii) follows. We now turn to prove~(i).
Since $\theta=\alpha+\beta$ is injective, $g_1g_2=0$. So we write $h=p-(g_1+g_2)$ and
claim that $h=0$. Write, for $i\in I, \ \tilde u_i=hu_i$ and $\tilde v_i=T(\tilde u_i)$
and note that Lemma~2.7 and Lemma~2.8 apply to $\tilde u_i,\tilde v_i$ in place of
$u_i,v_i$ (since $\{\tilde u_i\tilde u^*_i\}$ is an orthogonal family and
$\tilde u^*_i\tilde u_i=\tilde u^*_j\tilde u_j$ for all $i,j$). We also write
$\tilde r_i=\tilde v_i\tilde v^*_i, \ d_i=\tilde v^*_i\tilde v_i$.
Using Lemma~2.10 applied to $T^{-1}$, we find that $\tilde v_i=cv_i$ for some
projection $c\le v_iv^*_i=r_i$ and thus $\tilde r_i\in r_iNr_i$ and, by Lemma~2.11,
$\tilde r_i$ commutes with all $r_j$. In particular
$\tilde r_ir_0=r_0\tilde r_i, \ i\in I$.

Now $r_0\tilde v_i=\tilde v_i\tilde v_ir_0\tilde v_i\in \tilde r_iN\tilde d_i=
T(hu_iu^*_iMu^*_iu_i)$ and we can find $a=u_iu^*_iau^*_iu_i$ such that
$r_0\tilde v_i=T(ha)$. Hence
$$
\alpha(hau^*_i)=(r_i-r_0) T(ha)v^*_i=(r_i-r_0)r_0\tilde v_i\tilde v_i^*=0.
$$
But $\alpha$, restricted to $hMh$ is injective. Hence $hau^*_i=0$ and also $ha=0$ and
$r_0\tilde v_i=T(ha)=0$. But then $r_0\tilde r_i=0$ for all $i$; i.e.
$\tilde r_i\le r_i-r_0$ and it follows that $\tilde r_i\tilde r_j=0$ for all $i\ne j$.
Similar argument shows that $\tilde d_i\tilde d_j=0, \ i\ne j$. For a given $i\in I$
and $j\ne i$,
$$
\tilde r_i=\tilde r_i(I-\tilde r_j)=\tilde v_i\tilde d_j\tilde v_i^*=
\tilde v_i\tilde d_j\tilde d_i\tilde v_i^*=0.
$$
But then $\tilde v_i=0$ and, thus, $hu_i=0$ for all $i\in I$.
This shows that $h=0$ and we are done.

It is now left to prove (iii): $c(g_1)+c(g_2)=I$. But since $g_1+g_2=p$ and $c(p)=I$
and $g_i\in Z(pMp)$ it is obvious. \ \ \ \ $\blacksquare$
\medskip

Because of Lemma~2.4 it will suffice now, in order to prove Theorem~2.1 to restrict
our attention to the cases $c(g_1)=I$ and $c(g_2)=I$. Since the proof is similar in
these cases we now assume $c(g_1)=I$ (i.e. $g_2=c(g_2)=0$).

%                            LEMMA 2.20
\begin{lemma}
When $c(g_1)=I$ (with $g_1$ as in Lemma~2.19) we have, for all \linebreak
$a\in pMp$ and $x\in pM(I-p)$,
$$
T(ax)=\alpha(a)T(x).
$$
Moreover, we have now $r_0=0, \ \alpha$ is a $^*$-isomorphism of $pMp$ onto $rNr$ and
$T$ maps $pM(I-p)$ onto $rN(I-r)$ (i.e. $r=q$).
\end{lemma}

%                              PROOF
\noindent
{\bf Proof.}~~Now, that $g_2=0, \ \alpha$ is injective and, thus, a $^*$-isomorphism
onto $(r-r_0)N(r-r_0)$. But then we can repeat the argument of the proof of
Lemma~2.19(3)(ii) (with $p$ replacing $h$) to show that $r_0v_i=0$ for all $i$.
hence $r_0=0$ and $\alpha$ maps onto $rNr$. Now fix $x\in pM(I-p)$ with $x=u_iu^*_ix$.
Write $x=x_1+x_2$ with $x_1=u_iu^*_ixu^*_iu_i$ and $x_2=u_iu^*_ix(I-u^*_iu_i)$.
We have $T(x_1)\in T(u_iu^*_iMu^*_iu_i)=r_iNd_i$ (Lemma~2.9) and
$T(x_2)= T\bigl(u_iu^*_ix_2+x_2u^*_iu_i\bigr)=r_iT(x_2)+T(x_2)d_i$. However, for all
$j\ne i$, $x_2u^*_j=0$ (as $u^*_ju_j=u^*_iu_i$) and $u^*_jx_2=0$ (as $u_ju^*_ju_iu^*_i=0$);
hence (Lemma~2.9(3)) $T(x_2)v^*_j=0$ and also $T(x_2)d_j=T(x_2)v^*_jv_j=0$. But
$d_j=d_i$ (as $d_j(I-d_i)=v^*_jr_iv_j=v^*_jr_0v_j=0$) and we conclude that
$$
T(x_2)=r_iT(x_2)\subseteq rN.
$$
Hence $T(pM(I-p))\subseteq rN$. But $T(pM(I-p))=qN(I-q)$; hence $qN(I-q)\subseteq rN$
while $r\le q$. If $r\ne q$ then it follows from the fact that $c(I-q)=I$ that
$(q-r)N(I-q)\ne 0$ but this contradicts $qN(I-q)\subseteq rN$ and we get $q=r$.

We then conclude, from Corollary~2.18, that, given $i,j,k$ in $I, \ a\in u_iu^*_iMu_ju^*_j$
and $x\in u_ku^*_kM(I-p)$, we have
$$
T(ax)=\alpha(a)T(x).
$$
This equality then holds for finite sums of such $a,x$. Since $T$ is $\sigma$-weakly
continuous by \cite[Corollary~3.22]{Ho} and $\alpha$ is $\sigma$-weakly continuous
the equality holds for all $a\in pMp$ and $x\in pM(I-p)$. \ \ \ \ $\blacksquare$

%                              COROLLARY 2.21
\begin{corollary}
Assume $c(g_1)=I$. For all $x,y\in pM(I-p)$
$$
T(x)T(y)^*=\alpha(xy^*).
$$
\end{corollary}

%                              PROOF
\noindent
{\bf Proof.}~~This is \cite[Lemma~5.10]{MS} (which generalizes the result of
Lance~\cite[Theorem~3.5]{L2}). \ \ \ \ \ $\blacksquare$

%                                   PROPOSITION 2.22
\begin{proposition}
Assume $c(g_1)=I$ as above. Then there is a $^*$-isomorphism
$$
\gamma:(I-p)M(I-p)\longrightarrow (I-q)N(I-q)
$$
with
\begin{description}
\item(i) $\gamma$ is surjective
\item(ii) For $a\in(I-p)M(I-p)$ and $x\in pM(I-p)$,
$$
T(xa)=T(x)\gamma(a)
$$
\item(iii) For $x,y\in pM(I-p)$,
$$
T(x)^*T(y)=\gamma(x^*y)
$$
\end{description}
\end{proposition}

%                              PROOF
\noindent
{\bf Proof.}~~Suppose $N\subseteq B(H)$ (and the unit of $N$ is $I_H$) and write
$H_0=\overline{\rm span}\{T(x)^*h:x\in pM(I-p), \ h\in H\}$. Since
$T(pM(I-p))=qN(I-q)$ and $c(q)=I$, it follows that $H_0=(I-q)(H)$.

For $a\in(I-p)M(I-p)$ we define $\gamma(a)$ as an operator in $B(H_0)$ and assume that
it is defined to be zero on $H\ominus H_0$. We define
$$
\gamma(a)\left(\sum^n_{i=1} T(x_i)^*h_i\right)=\sum^n_{i=1} T(x_ia^*)^*h_i
$$
where $x_i\in pM(I-p), \ h_i\in H.$

Note the following
\be
\left<\sum T(x_ia^*)^*h_i,\sum T(x_ia^*)^*h_i\right>&=&
\sum_{i,j}\left<T(x_ja^*)T(x_ia^*)^*h_i,h_j\right>\nonumber\\
&=&\sum_{i,j}\left<\alpha(x_ja^*ax^*_i)h_i,h_j\right>.\nonumber
\ee
Since the matrix $(x_ja^*ax^*_i)\in M_n(pMp)$ is majorized by the matrix
$\|a\|^2(x_jx^*_i)$ and $\alpha$ is a $^*$-isomorphism,
$$
\bigl(\alpha(x_ja^*ax^*_i)\bigr)\le \|a\|^2\bigl(\alpha(x_jx^*_i)\bigr).
$$
Hence
\be
\left\|\sum T(x_ia^*)^*h_i\right\|^2&\le&
\|a\|^2\sum\left<\alpha(x_jx^*_i)h_i,h_j\right>\nonumber\\
&=&\left\|\sum T(x_i)^*h_i \right\|^2\|a\|^2.\nonumber
\ee
It follows that $\gamma(a)$ is well defined and can be extended to an
operator in $B(H)$ with $\|\gamma(a)\|\le \|a\|$. For
$x_i\in pM(I-p), \ h_i\in H, \ i=1,2$, we have
\be
&&\left<\gamma(a)T(x_i)^*h_1, \ T(x_2)^*h_2\right>=
\left<T(x_1a^*)^*h_1, \ T(x_2)^*h_2\right>\nonumber\\
&&=\left<\alpha(x_2ax^*_1)h_1,h_2\right>=\left<T(x_1)^*h_1,T(x_2a)^*h_2\right>\nonumber\\
&&=\left<T(x_1)^*h_1,\gamma(a^*)T(x_2)^*h_2\right>\nonumber
\ee
Hence $\gamma(a^*)=\gamma(a)^*$. It is easy to check that $\gamma$ is multiplicative
and injective.
Now suppose $b\in N'$ then for every $x\in pM(I-p), \ h\in H$,
\be
\gamma(a)bT(x)^*h=\gamma(a)T(x)^*bh=T(xa^*)^*bh=bT(xa^*)^*h=b\gamma(a)T(x)^*h.\nonumber
\ee
Hence $\gamma(a)\in N$. Since $\gamma(a)$ is zero on
$H\ominus H_0, \ \gamma(a)\in (I-q)N(I-q)$.

It follows from the definition that, for $a\in(I-p)M(I-p)$ and $x\in pM(I-p), \
\gamma(a)T(x)^*=T(xa^*)^*$; hence
$$
T(x)\gamma(a)=T(xa).
$$
This proves part~(ii).

Now choose $z,t\in qN(I-q)$ and write $a=T^{-1}(z)^*T^{-1}(t)\in (I-p)M(I-p)$.
Compute, for $x\in pM(I-p)$ and $h\in H$,
\be
\gamma(a)T(x)^*h&=&T\bigl(xT^{-1}(t)^*T^{-1}(z) \bigr)^*h=[\alpha(xT^{-1}(t)^*)z]^*h\nonumber\\
                &=&[T(x)t^*z]^*h=z^*tT(x)^*h=z^*tT(x)^*h.\nonumber
\ee
Hence $z^*t\in\gamma\bigl((I-p)M(I-p)\bigr)$. Since products of this form generate
$(I-q)N(I-q)$ as a von~Neumann algebra and the image of $\gamma$ is a von~Neumann
algebra, $\gamma$ is surjective. This proves (i). We can now apply \cite[Lemma~5.10]{MS}
to get~(iii). \ \ \ \ \ $\blacksquare$
\medskip

%                                 REMARK
\noindent
{\bf Remark.}~~Note that $\gamma$ is in fact equivalent to the representation on the
internal tensor product $(I-p)Mp\otimes_\alpha H$.

%                              COROLLARY 2.23
\begin{corollary}

Assume $c(g_1)=I$. Then there is a $^*$-isomorphism
$$
\Psi:M\longrightarrow N
$$
such that, for $a\in M$,
\begin{description}
\item(i)~~~$\Psi(pap)=\alpha(pap)\in qNq$
\item(ii)~~$\Psi(pa(I-p))=T(pa(I-p))\in qN(I-q)$
\item(iii)~$\Psi((I-p)a(I-p))=\gamma((I-p)a(I-p))\in (I-q)N(I-q)$
\item(iv)~~$\Psi((I-p)ap)=T(pa^*(I-p))^*\in(I-q)Nq.$
\end{description}
\end{corollary}

%                              PROOF
\noindent
{\bf Proof.}~~The equations (i)-(iv) define $\Psi$ and the properties of $\alpha$ and
$\gamma$ (see Lemma~2.20, Corollary~2.21 and Proposition~2.22) show that $\Psi$ is
indeed a $^*$-isomorphism of $M$ onto $N$. \ \ \ \ \ $\blacksquare$
\medskip

%              PROOF of Theorem 2.1

\noindent
{\bf Proof of Theorem~2.1.}~~It was shown in Lemma~2.4 that, to prove Theorem~2.1, it
will suffice to write the algebra as a direct sum (using central projections) of
algebras for which Theorem~2.1 holds. In Lemma~2.5 we saw that we can assume that the
algebra $M$ and the projections $p$ and $I-p$ satisfy one of the conditions ((1), (2)
or (3)) stated in that lemma. If condition~(1) is satisfied then Theorem~2.1 follows
from Lemma~2.6. So we can assume that either condition~(2) or condition~(3) is
satisfied. Condition~(3) is, in fact, condition~(2) for $p,I-p$ interchanged. It
will suffice, therefore, to assume condition~(2). We then find, in Lemma~2.19, two
central projections, $c(g_1)$ and $c(g_2)$, with $c(g_1)+c(g_2)=I$ and (again, by
referring to Lemma~2.4) we can assume that either $c(g_1)=I$ or $c(g_2)=I$. For
the first case the theorem is proved in Corollary~2.23. In this case we get, in fact,
that $e_2=0$ and the map, extending $T$, is a $^*$-isomorphism. The proof of the other
case, when $c(g_2)=I$, is similar and is omitted. In that case the map turns out to
be a $^*$-antiisomorphism. \ \ \ \ \ $\blacksquare$

%                                 REMARK 2.24

\begin{remark}
The map $\Phi+\Psi$ of Theorem~2.1 is an isometry of $M$ onto $N$ that extends $T$ and maps
$pMp+(I-p)M(I-p)$ onto $qNq+(I-q)N(I-q)$ and $pM(I-p)$ onto $qN(I-q)$. We write
$\Lambda$ for $\Phi+\Psi$.
\end{remark}

Recall from the introduction that given a right selfdual Hilbert $C^*$-module $X$ over
a von~Neumann algebra $A$, we can form the von~Neumann linking algebra which can be
written
$$
{\cal L}_w(X)=
\begin{pmatrix}
\BB\, (X) & X\\
\bar X  & A
\end{pmatrix}
$$
where $\BB\,(X)$ is the algebra of all bounded, adjointable $A$-linear maps on $X$ and
$\bar X$ is the conjugate module (which is a left Hilbert $C^*$-module over $A$).
It is known that this algebra is indeed a von~Neumann algebra. We assume that our
$C^*$-modules are full and this implies that we can write $X$ as $p{\cal L}_w(X)(I-p)$
for a projection $p$ with $c(p)=c(I-p)=I$. The following corollary then follows
immediately from the theorem.

%                              COROLLARY 2.25
\begin{corollary}
If $X$ and $Y$ are right selfdual $C^*$-modules over (possibly different) von~Neumann
algebras $A,B$ and if ${\cal L}_w(X)$ and ${\cal L}_w(Y)$ are the von Neumann linking
algebras of $X$ and $Y$ respectively then every linear surjective isometry %T:X\to Y$
extends to a linear surjective isometry $\Lambda:{\cal L}_w(X)\to {\cal L}_w(Y)$.

Moreover, there is a central projection $z\in{\cal L}_w(X)$ such that if we write
$\Psi=\Lambda|z{\cal L}_w(X)$ and $\Phi=\Lambda|(I-z){\cal L}_w(X)$ then
\medskip

\noindent
(1)~~~$\Psi$ is a $^*$-isomorphism onto $\Lambda(z){\cal L}_w(Y)$. It defines a
$^*$-isomorphism \linebreak
$\Psi_{11}:z\BB(X)\to \Lambda(z)\BB(Y)$ and a $^*$-isomorphism $\Psi_{22}:zA\to \Lambda(z)B$
such that, for $L\in\BB(X)z, \ a\in Az$ and $x\in zX, \ y\in zX$,
\be
&&T(Lxa)=\Psi_{11}(L)T(x)\Psi_{22}(a)\nonumber\\
&&T(x)T(y)^*=\Psi_{11}(x\otimes y^*)\nonumber\\
&&T(x)^*T(y)=\Psi_{22}(\left<x,y\right>_A)\nonumber
\ee
and
\medskip

\noindent
(2)~~~$\Phi$ is a $^*$-antiisomorphism onto $\Lambda(I-z){\cal L}_w(Y)$. It defines
$^*$-antiisomorphisms $\Phi_{11}:(I-z)\BB\,(X)\to\Lambda(I-z)B$ and
$\Phi_{22}:(I-z)A\to \Lambda(I-z)\BB\,(Y)$ such that, for $L\in \BB\,(X)(I-z), \
x,y,\in (I-z)X, \ a\in(I-z)A,$
\be
&&T(Lxa)=\Phi_{22}(a)T(x)\Phi_{11}(L)\nonumber\\
&&T(x)T(y)^*=\Phi_{22}(\left<x,y\right>)\nonumber\\
&&T(x)^*T(y)=\Phi_{11}(y\otimes x^*).\nonumber \ \ \ \ \ \blacksquare
\ee

\end{corollary}

%                              COROLLARY 2.26
\begin{corollary}
If $X,Y,A,B$ are as in Corollary~2.24 but we assume also that $A$ is a factor then
the surjective linear isometry $T:X\to Y$ can be extended to a map
$$
\Lambda:{\cal L}_w(X)\to{\cal L}_w(Y)
$$
which is either a $^*$-isomorphism of or a $^*$-antiisomorphism.

Moreover, if $X^t$ is any $C^*$-module that is isometric to $X$ and such that the
induced isometry from ${\cal L}_w(X^t)$ to ${\cal L}_w(X)$ is a $^*$-antiisomorphism,
then $Y$ is completely isometric  to either $X$ or $X^t$. \ \ \ \ \ $\blacksquare$
\end{corollary}

\section{Isometries of Hilbert $C^*$-modules}

Now let $X$ and $Y$ be right (full) Hilbert $C^*$-modules over the $C^*$-algebras $A$ and
$B$ respectively and let $T:X\to Y$ be a surjective linear isometry. Write ${\cal L}(X)$
(and ${\cal L}(Y)$) for the linking algebra of $X$ (and $Y$); i.e.
$$
{\cal L}(X)=
\begin{pmatrix}
\KK(X) &X\\
\bar X & A
\end{pmatrix};\
{\cal L}(Y)=
\begin{pmatrix}
\KK(Y) & Y\\
\bar Y & B
\end{pmatrix}
$$
where $\bar X$ and $\bar Y$ are the conjugate modules. $\bar X$ is a left Hilbert
$C^*$-module over $A$ and $\bar Y$ is over $B$. (See Section~1 for the definitions).
It is known that ${\cal L}(X)$ and ${\cal L}(Y)$ are $C^*$-algebras. Now the isometry
$T:X\to Y$ induces an isometry $T^{**}:X^{**}\to Y^{**}$ (where $X^{**}$ is the second
dual of $X$) extending $T$. Write
$$
M={\cal L}(X)^{**}, \ N={\cal L}(Y)^{**}
$$
Then $M$ (respectively $N$) is a von~Neumann algebra. In fact it can be identified
with the universal enveloping  algebra of ${\cal L}(X)$ (respectively ${\cal L}(Y)$).

To continue we note the following.

%                             LEMMA 3.1
\begin{lemma}
There is a projection $p\in {\cal L}(X)^{**}$ such that $p{\cal L}(X)^{**}(I-p)$ is
the $w^*$-closure of $X\subseteq{\cal L}(X)\subseteq{\cal L}(X)^{**}$. Similarly,
the $w^*$-closure of $Y$ in ${\cal L}(Y)^{**}$ is $q{\cal L}(Y)^{**}(I-p)$ for some
projection $q\in{\cal L}(Y)^{**}$.

Moreover, $c(p)=c(I-p)=I$ and similarly for $q$.
\end{lemma}

%                             PROOF
\noindent
{\bf Proof.}~~We can identify ${\cal L}(X)^{**}$ with the $\sigma$-weak closure of
$\pi_u({\cal L}(X))$ where $(\pi_u,H_u)$ is the universal representation of
${\cal L}(X)$. The $w^*$-topology of ${\cal L}(X)^{**}$ is then the $\sigma$-weak
topology of $\pi_u({\cal L}(X))$. Write
$$
{\cal I}_1=
\begin{pmatrix}
\KK(X) & X\\
0 & 0
\end{pmatrix}
\subseteq{\cal L}(X); \ \  {\cal I}_2=
\begin{pmatrix}
0 & 0\\
\bar X & A
\end{pmatrix}
\subseteq {\cal L}(X).
$$
Write $M$ for the $\sigma$-weak closure of $\pi_u({\cal L}(X))$. ${\cal I}_1$ and
${\cal I}_2$ are right ideals in ${\cal L}(X)$ and the $\sigma$-weak closures of
${\cal I}_1$ and ${\cal I}_2$ are of the form $p_1M$ and $p_2M$ respectively. In fact,
$p_1=\vee\{r(y):y\in\pi_u({\cal I}_1)\}$ and $p_2=\vee\{r(z):z\in\pi_u({\cal I}_2)\}$
where $r(y)$ is the range projection (in $M$) of $y$. The $\sigma$-weak closure of
$\pi_u(X)$ is then $p_1M\cap(p_2M)^*=p_1Mp_2$. But it is clear that $p_1+p_2=I$ and,
thus, writing $p=p_1$ we find that the $\sigma$-weak closure of $\pi_u(X)$ is $pM(I-p)$.

A similar argument works for $Y$.

Finally note that $MpM$ contains $\pi_u({\cal L}(X){\cal I}_1)=\pi_u({\cal L}(X))$;
hence is $\sigma$-weakly dense in $M$ and it follows that $c(p)=I$. The argument for
$I-p(=p_2)$ is similar. \ \ \ \ \ $\blacksquare$
\medskip

Since it follows from Banach space theory that $X^{**}$ is isometrically isomorphic
to the $w^*$-closure of $X$ in ${\cal L}(X)^{**}$ we conclude that $T^{**}$ induces
a surjective linear isometry
$$
S:p{\cal L}(X)^{**}(I-p)\to q{\cal L}(Y)^{**}(I-q)
$$
and $S$ extends $T$ (when we view $X,Y$ as subspaces of ${\cal L}(X)^{**}$ and
${\cal L}(Y)^{**}$ respectively). In particular, $S$ maps $X$ onto $Y$.

Applying Theorem~2.1 to $S$ we obtain the following.

%                              THEOREM 3.2

\begin{theorem}
Let $X$ be a right full Hilbert $C^*$-module over the $C^*$-algebra $A$ and $Y$ be a
right full Hilbert $C^*$-module over $B$. Let $T:X\to Y$ be a surjective linear
isometry. Then
\begin{description}
\item(1) There is a surjective linear isometry
$$
\Lambda_0:{\cal L}(X)\to {\cal L}(Y)
$$
extending $T$ and mapping $\KK(X)\oplus A$ onto $\KK(Y)\oplus B$.

\item(2) There is a projection $f$ in the center of ${\cal L}(Y)^{**}$ such that
the map
$$
\Psi_0(a)=\Lambda_0(a)f \ \ \ a\in {\cal L}(X)
$$
is a $^*$-homomorphism and the map
$$
\Phi_0(a)=\Lambda_0(a)(I-f) \ \ \ a\in {\cal L}(X)
$$
is a $^*$-antihomomorphism.
\end{description}

\end{theorem}

%                             PROOF
\noindent
{\bf Proof.}~~Part~(2) follows from Theorem~2.1 but in fact it is known whenever
$\Lambda_0$ is an isometry of $C^*$-algebras (see~\cite{K} and \cite[p.\,188]{T}).
For part~(1) we need only to notice that the isometry $\Lambda=\Psi+\Phi$, given by
Theorem~2.1 and Remark~2.24 (mapping ${\cal L}(X)^{**}$ onto ${\cal L}(Y)^{**}$)
maps ${\cal L}(X)$ onto ${\cal L}(Y)$. But we know from Corollary~2.25 that $\Psi$
maps $\KK(X)$ onto $\KK(Y)$ since
$$
\Psi(x\otimes y^*)=S(x)\otimes S(y)^*=T(x)\otimes T(y)^*\in\KK(Y)
$$
and similarly $\Psi$ maps $A$ onto $B$, $X$ onto $Y$ and $\bar X$ onto $\bar Y$.
The statements for $\Phi$ are similar (although here
$\Phi(\KK(X))=B, \ \Phi(A)=\KK(Y))$. \ \ \ \ \ $\blacksquare$

%                                 REMARK 3.3
\begin{remark}
If we have $X,Y$ and $T$ as in Theorem~3.2 and we let $\Psi_{11}$ define the
restriction of $\Psi_0$ to $\KK(X)$ and $\Psi_{22}$ be the restriction of $\Psi_0$ to
$A$, then it follows from the properties of $\Psi_0$ that both $\Psi_{11}$ and $\Psi_{22}$
are $^*$-homomorphisms and
$$
T(Kxa)f=\Psi_{11}(K)T(x)\Psi_{22}(a) \ \ \ a\in A, \ K\in\KK(X), \ x\in X.
$$
Similarly, we get $^*$-antihomomorphisms $\Phi_{11}$ and $\Phi_{22}$ with
$$
T(Kxa)(I-f)=\Phi_{22}(a)T(x)\Phi_{11}(K) \ \ \ a\in A, \ K\in\KK(X), \ x\in X.
$$
\end{remark}

%                              COROLLARY 3.3
\begin{corollary}

Suppose $X$ and $Y$ are as in Theorem~3.2 and $T:X\to Y$ is a 2-isometry (i.e. the map
$I\otimes T$ that maps $M_2\otimes X$ onto $M_2\otimes Y$ is an isometry). Then the map
$\Lambda_0$ of Theorem~3.2 is a $^*$-isomorphism (i.e. $\Lambda_0=\Psi_0, \ f=I$) and
if we let $\Psi_{11}$ and $\Psi_{22}$ be the maps induced by $\Lambda_0=\Psi_0$ on
$\KK(X)$ and on $A$ respectively (so that $\Psi_{11}:\KK(X)\to \KK(Y)$ and
$\Psi_{22}:A\to B$) then $\Psi_{11}$ and $\Psi_{22}$ are $^*$-isomorphisms and
\begin{description}
\item(i)~~~$T(Lxa)=\Psi_{11}(L)T(x)\Psi_{22}(a), \ \ \ L\in\KK(X), \ a\in A, \ x\in X$
\item(ii)~~$T(x)^*T(y)=\Psi_{22}(\left<x,y \right>) \ \ \ x,y\in X$
\item(iii)~$T(x)T(y)^*=\Psi_{11}(x\otimes y^*) \ \ \ x,y\in X$.
\end{description}
\end{corollary}

%                             PROOF
\noindent
{\bf Proof.}~~Write $T_2$ for the isometry $T_2:M_2\otimes X\to M_2\otimes Y$. Fix
$x_0,y_0,z_0\in X$ and write
$$
x=e_{12}\otimes x_0, \ y=e_{11}\otimes y_0, \ z=e_{21}\otimes z_0 \ \ {\rm in} \ \ M_2\otimes X.
$$
Since $T_2$ satisfies Lemma~2.2,
$$
T_2(xy^*z+zy^*x)=T_2(x)T_2(y)^*T_2(z)+T_2(z)T_2(y)^*T_2(x).
$$
But $xy^*z=0$ and $T_2(x)T_2(y)^*T_2(z)=0$. This implies that
$T_2\bigl(e_{22}\otimes z_0y^*_0x_0 \bigr)=T_2(zy^*x)=T_2(z)T_2(y)^*T_2(x)=
e_{22}\otimes T(z_0)T(y_0)^*T(x_0)$. Hence $T(z_0y^*_0x_0)=T(z_0)T(y_0)^*T(x_0)$ for all
$x_0,y_0,z_0$ in $X$.

We again write $X^{**}$ as $pM(I-p) \ (M={\cal L}(X)^{**})$ and
$Y^{**}$ as $qN(I-q) \ (N={\cal L}(Y)^{**})$. We get an isometry $S:pM(I-p)\to qN(I-q)$
extending $T$. Clearly this extension still satisfies
$$
S(zy^*x)=S(z)S(y)^*S(x) \ \ \ x,y,z\in pM(I-p).
$$
We can now apply the results of the previous section to $S$.

As in Lemma~2.5, $M$ can be decomposed into a direct sum of von~Neumann algebras
each satisfying one of the conditions stated in Lemma~2.5. If condition~(1) holds
then, using the proof of lemma~2.6 we see that the map induced on the linking
algebra is a $^*$-isomorphism. Suppose now that condition~(2) holds. Let $\{u_i\}$
be as in this condition and write $v_i=S(u_i), \ r_i=S(u_i)S(u_i)^*$ and
$d_i=S(u_i)^*S(u_i)$. We have, for $i\ne j$, $$
0=S(u_iu^*_iu_j)=S(u_i)S(u_i)^*S(u_j). $$

Hence $r_ir_j=0$. It follows form Lemma~2.7 that $d_i=d_j$ for all $i,j$. It then
follows that the map $\beta$, defined in the discussion following Lemma~2.13, vanishes;
i.e. $\theta=\alpha$. Hence the map induced on the linking algebra $M$ is a $^*$-isomorphism.
A similar argument works if condition~(3) (of Lemma~2.5) holds.

Statements (i)-(iii) now follow.  \ \ \ \ \ $\blacksquare$

%                            REMARK 3.5
\begin{remark}
Corollary~3.4 shows that a (surjective) 2-isometry from one $C^*$-Hilbert module to
another is necessarily completely isometric. I do not know to what extent it holds
for larger classes of operator spaces. Argument as above (combined with
Corollary~2.10 of~\cite{AS}) can be used to show that unital 2-isometries of
operator algebras are multiplicative and for some classes of operator algebras this
would imply that they are complete isometries.
\end{remark}

Given a Hilbert space $H$ there is more than one way to represent it as an operator
space; i.e. one can find different operator spaces that are all isometric to $H$ but
they are pairwise non-completely-isometric.

One representation of $H$ as an operator space is when you fix an orthonormal basis
$\{e_i\}$ for $H$ and consider the space of all bounded operators in $B(H)$ whose
matrix with respect to this basis has non zero entries only in the first column. This
subspace of $B(H)$ is isometric to $H$ and is called the Hilbert column space, written
$H^c$. One can write $H^c=B(\CC\, e_1,H)$.

Replacing the word ``column'' by the word ``row'' we  get the Hilbert row space $H^r$
(or $B(H,\CC\, e_1)$). It is known that these two operator spaces are isometric but not
completely isometric. Both operator spaces have a natural Hilbert $C^*$-module
structure. $H^c$ is a $C^*$-module over the algebra $\CC$ and $H^r$ is a $C^*$-module
over the algebra $K(H)$, the compact operators on $H$.

But in addition to $H^c$ and $H^r$ there are many other different (i.e. not completely
isometric) representations of $H$ as an operator space (see~\cite{Pi}). The following
corollary shows that none of these is a Hilbert $C^*$-module.

%                              COROLLARY 3.5
\begin{corollary}
Let $H$ be a fixed Hilbert space and $X$ be a Hilbert $C^*$-module over a
$C^*$-algebra $A$ that is isometric (as a Banach space) to $H$. Then either $A$ is
isomorphic to $\CC$ and $X$ is completely isometrically isomorphic to $H^c$ or $A$ is
$^*$-isomorphic to $K(H)$ and $X$ is completely isometrically isomorphic to $H^r$.
\end{corollary}

%                             PROOF
\noindent
{\bf Proof.}~~Assume, for simplicity, that $H$ is infinite dimensional. Write
$T:H^c\to X$ for the linear surjective isometry of $H^c$ onto $X$. From Theorem~3.2 it
follows that we can extend $T$ to an isometry
$$
\Lambda_0:{\cal L}(H^c)=
\begin{pmatrix}
K(H) & H^c\\
\overline{H^c} & \CC
\end{pmatrix}
\longrightarrow {\cal L}(X)=
\begin{pmatrix}
\KK(X) & X\\
\bar X & A
\end{pmatrix}.
$$
But ${\cal L}(H^c)$ is $^*$-isomorphic to $K(H)$ and thus ${\cal L}(H^c)^{**}$ is a
factor. Hence $\Lambda_0$ is either a $^*$-isomorphism or a $^*$-antiisomorphism.
In the former case $\Lambda_0$ maps $\CC$ onto $A$ and $H^c$ onto $X$ and is a
complete isometry. In the latter case we consider the map $\tau$
$$
\tau:
\begin{pmatrix}
\CC & H^r\\
\overline{H^r} & K(H)
\end{pmatrix}
\longrightarrow
\begin{pmatrix}
K(H) & H^c\\
\overline{H^c} & \CC
\end{pmatrix}
$$
defined by
$$
\tau
\begin{pmatrix}
\lambda & y\\
\bar z & K
\end{pmatrix}
=
\begin{pmatrix}
K^t & y^t\\
\overline{z^t} & \lambda
\end{pmatrix}
\ \ \ \ y,z\in H^r, \ K\in K(H), \ \lambda\in\CC.
$$
(where $K^t$ is the transpose of $K$ and $y^t$ is the transpose of $y$). Then
$\tau$ is a $^*$-antiisomorphism and $\Lambda_0\circ\tau$ is then a $^*$-isomorphism
that maps, completely isometrically, $K(H)$ onto $A$ and $H^r$ onto $X$. \ \ \ \
$\blacksquare$

%++++++++++++++++++++++++++++++++++++++++++++++++++++++++++++++++++++++++

%                           REFERENCES

%\baselineskip 14pt

\end{document}